\documentclass[11pt,twoside]{report}

\usepackage{a4}
\usepackage{amsmath}
\usepackage{amssymb}
\usepackage{diagrams}

\DeclareMathAlphabet{\scr}{U}{eus}{m}{n}

\newcommand\Z{{\mathbb Z}}
\newcommand\Q{{\mathbb Q}}

\newcommand\PP{{\mathbb P}}
\newcommand\A{{\mathbb A}}
\newcommand\G{{\mathbb G}}

\newcommand\cF{{\cal F}}

\newcommand\cI{{\cal I}}
\newcommand\cO{{\cal O}}

\newcommand\sL{{\cal L}} % sollte "Schreibschrift" sein

\newcommand\Mloc{{\rm\bf M}^{\rm loc}}
\newcommand\Mlocs{\overline{\rm\bf M}^{\rm loc}} % spez. Faser
\newcommand\Mloct{\widetilde{\rm \bf M}^{\rm loc}}
\newcommand\ii{{\rm\bf i}}

\newcommand\xt{[\hspace{-.12em}[ t ]\hspace{-.12em} ] } 
\newcommand\xT{(\!( t )\!)}

\newcommand\tens\otimes

\newcommand\lto{\longrightarrow}
\newcommand\id{\mathop{\rm id}\nolimits}
\renewcommand{\div}{{\rm div}}
\renewcommand{\hom}{\mathop{\rm Hom}\nolimits}

\newcommand{\Char}{\mathop{\rm char}\nolimits}
\newcommand{\Spec}{\mathop{\rm Spec}}

\newcommand{\Grass}{\mathop{\rm Grass}\nolimits}
\newcommand{\Flag}{\mathop{\rm Flag}\nolimits}
\newcommand{\Mat}{\mathop{\rm Mat}\nolimits}
\newcommand{\Res}{\mathop{\rm Res}\nolimits}
\renewcommand{\det}{\mathop{\rm det}\nolimits}
\newcommand{\limind}{{{\mathop{\rm lim}}_{\rightarrow}}}
\newcommand{\pr}{\mathop{\rm pr}}

\newcommand\qed{\hfill$\square$}

\newtheorem{thm}{Theorem}[section]
\newtheorem{stz}[thm]{Proposition}
\newtheorem{lem}[thm]{Lemma}
\newtheorem{Def}[thm]{Definition}
\newtheorem{kor}[thm]{Corollary}
\newtheorem{bed}[thm]{Condition}

\begin{document}

\begin{titlepage}
\hspace*{1cm}
\vskip3cm

\begin{center}

{\bf \huge On the flatness of models } 

{\bf \huge of certain Shimura varieties} 

{\bf \huge of PEL-type }

\vskip5cm

{\bf\large Ulrich G\"ortz}\\[9cm]

{\bf\large December 1999}
\end{center}
\end{titlepage}

\setcounter{tocdepth}{1}

\tableofcontents

\chapter{Introduction}

In the arithmetic theory of Shimura varieties, it is of interest
to have a model of the Shimura variety 
over the ring of integers $O_{E}$,
where $E$ is the completion of the reflex field at some
place lying over a prime $p$. For a Shimura variety of PEL-type,
which is a moduli
space of abelian varieties with certain additional
structure, it is a natural idea to define such a model 
by posing the moduli problem over $O_{E}$. In the case of 
a hyperspecial level structure, one gets a smooth model
as was shown by Kottwitz.
% \cite{K}.
We are interested in the case of parahoric level structures,
where one cannot expect to get a smooth model.

The problem reduces in some sense to the level of $p$-divisible
groups.
In their book \cite{RZ}, Rapoport and Zink
investigate (formal) moduli schemes of $p$-divisible groups
and define such models.
Unfortunately,
very little is known in general about the structure of these models.
In most cases one does not even know if they are flat over $O_{E}$ ---
which is certainly a condition a reasonable model should 
satisfy. 

A different approach to these problems is to look for another model
which is semi-stable or at least has toroidal singularities.
There is an interesting proposal of Genestier of
a semi-stable model in the case of the symplectic group
(cf. \cite{G}), which works in low dimensions. 
Recently Faltings announced some results about a model with 
toroidal singularities, which again works in low
dimensions (cf. \cite{F2}).
On the other hand, a drawback of this approach
is that the new model does not have an easy description
as a moduli space of $p$-divisible groups.

To examine local properties such as flatness, it is
useful to work with the so-called {\em local model},  which
locally for the \'etale topology around each point
of the special fibre coincides with the
corresponding moduli scheme of $p$-divisible groups,
but which can be defined in terms of linear algebra
and is thus much easier to handle (cf. \cite{RZ}).

In this article, we will deal with the flatness conjecture
in a special case which will be explained now.
We use the same notation as in \cite{RZ}.

Let $F/\Q_p$ be a finite {\em unramified} extension,
let $B=F$, $V=F^n$. The algebraic group associated
to these data is $G = \Res_{F/\Q_p} GL_F(V)$.

Let $L$ be an algebraically closed field of characteristic $p$,
and denote by $K_0$ the quotient field of the Witt ring $W(L)$.
As $F$ is unramified over $\Q_p$, with the notation of \cite{RZ}
we have $K_0=K$.
 
Furthermore let $\mu: \G_{m,K} \lto G_K$ be a 1-parameter subgroup,
such that the weight decomposition of $V \tens_{\Q_p} K$ contains
only the weights 0 and 1: 
$$V\tens_{\Q_p} K = V_0 \oplus V_1.$$
Denote by $E$ the field of definition of the conjugacy
class of $\mu$.

Finally, let $\sL$ be a periodic lattice chain
in $V$.

We will call data of this type {\em unramified data of (EL) type}. 
Note that this is a more general notion than is used in 
\cite{RZ}, 3.82. There, the lattice chain consists only of
multiples of 
one lattice, and the resulting local model is smooth.
Here, 'unramified' relates only to the field extension.

Given these data, Rapoport and Zink define a local model 
(see \cite{RZ} 3.27 or section \ref{unram_ext}),
which is a projective $O_E$-scheme. We will give the definition
in the case $F=\Q_p$ below.

These local models are related to Shimura varieties of 
PEL-type which belong to unitary groups
that split over an unramified extension of $\Q_p$.

Our main result is the following theorem which confirms the
conjecture of Rapoport and Zink in this case (see theorem
\ref{general_case}).

{\bf Main Theorem.} {\em The local model associated 
to an unramified (EL)-datum is flat over $O_E$, 
and its special fibre is reduced.
The irreducible components of the special fibre
are normal with rational singularities,
so in particular are Cohen-Macaulay.}

It is essential that we consider only unramified extensions
 $F/\Q_p$. In fact, in the ramified case
 the flatness conjecture has to be refined
as is shown by the results of Pappas \cite{P}.

It is an interesting question if the special fibre as a whole
has Cohen-Macaulay singularities. In view of the flatness this is
equivalent to the local model being Cohen-Macaulay. The 
second remark after proposition \ref{equ_Utau} shows that 
this would follow if
one could prove that the affine scheme (over some field)
 defined by the equations
$$ B_{m-1}B_{m-2}\cdots B_0 =  B_{m-2}\cdots B_0 B_{m-1} = \cdots 
   = B_0 B_{m-1} \cdots B_1 = 0, $$
where the $B_i$ are generic $k \times k$-matrices, 
is Cohen-Macaulay and has the 'right' dimension.

Of course, the flatness question presents itself also for other
groups, in particular for the symplectic group $Sp_{2n}$. 
% In terms of the local model,
% this means that the $\Lambda_i$ are additionally equipped with an
% alternating form and that the $\cF_i$ have to be 'compatible'
% with this form.
In this situation flatness has been verified in special cases
by Deligne and Pappas \cite{DP}, de Jong  \cite{dJ},
and Chai and Norman \cite{CN}. 
Their proofs rely heavily on
very explicit calculations with the equations. 
One of the theories
involved is the theory of algebras with straightening law. 
It allows one to show in some cases that the special fibre
of the local model is reduced, or even that it is Cohen-Macaulay.
The disadvantage of this method is that only cases
where the lattice chain is
small, i. e. does not consist of too many lattices, can be handled. 

To give an idea of the proof of the main theorem, let us 
first give the definition of the standard local model,
where in terms of the (EL) datum we have $F=\Q_p$
(compare section \ref{def_Mloc}). In fact,
it is easy to see that the local model associated to an 
unramified (EL)-datum is isomorphic
after unramified base change to a product of standard local 
models. 

Let $O$ be a complete discrete valuation ring with
perfect residue class field. Let $\pi$ be a
uniformizer of $O$ and denote the field of fractions of $O$
by $K$. Let $k$ be an algebraic closure of the residue class field
of $O$.
Fix integers $0< r < n$.
Let $e_1, \dots, e_n$
be the canonical basis of $K^n$.
Let $\Lambda_i$, $0 \le i \le n-1$, be the free $O$-module of rank 
$n$ with basis 
 $e^i_1 := \pi^{-1} e_1, \dots,e^i_i := \pi^{-1} e_{i},e^i_{i+1} :=e_{i+1}, 
   \dots, e^i_n := e_n$.
This yields a complete lattice chain
\begin{diagram}
\cdots & \rTo & \Lambda_0 & \rTo & \Lambda_1 
& \rTo & \cdots & \rTo & \Lambda_{n-1} & \rTo & \pi^{-1}\Lambda_0 & \cdots
\end{diagram}
Choose $I = \{i_0 < i_1 < \cdots < i_{m-1} \} \subseteq \{0, \dots, n-1\}$.

Then the standard local model $\Mloc_I$ is the $O$-scheme
that represents the following functor. For an $O$-scheme
 $S$, the $S$-valued points of $\Mloc_I$ are
the isomorphism classes
of commutative diagrams 
\begin{diagram}
\Lambda_{i_0, S} & \rTo & \Lambda_{i_1,S} & \rTo & \cdots 
& \rTo & \Lambda_{i_{m-1},S} & \rTo^\pi & \Lambda_{i_0,S} \\
\uInto              &      & \uInto        &      &   &&\uInto & &\uInto   \\
\cF_0     & \rTo & \cF_1  & \rTo & \cdots & \rTo & \cF_{m-1}
  &\rTo & \cF_0
\end{diagram}
where $\Lambda_{i,S}$ is $\Lambda_i \tens_{O} \cO_S$, and where 
the $\cF_\kappa$ are locally free $\cO_S$-submodules of rank $r$ which
Zariski-locally on $S$ are direct summands of $\Lambda_{i_\kappa,S}$. 
We write $\Mloc := \Mloc_{\{0, \dots, n-1 \}}$.

We see that this functor is representable by a closed
subscheme of a product of Grassmannians. The generic
fibre is a Grassmannian itself since all the maps
 $\Lambda_i \lto \Lambda_j$ are isomorphisms after 
tensoring with $K$.

The case where $r= n-1$ (the so-called Drinfeld case)
is particularly simple: then $\Mloc$ has semi-stable reduction, 
and $\Mloc$ is
obviously flat. 
On the other hand, if $I$ is small, the local model
is less complicated than for large $I$. For example if
 $I$ has only one element, then $\Mloc_I$ is simply a Grassmannian 
over $O$, so it is even smooth.
In this work we will show that $\Mloc_I$ is flat over $O$
for general $n$, $r$ and $I$.

The most difficult part is to show that the special 
fibre of the local model is reduced.
In positive characteristic, this question
can be reduced to a question on local models for 
small $m$ by embedding
the special fibre of the 
local model into the affine flag variety and using
the technique of Frobenius splitting. 
Let us make this a little more precise.

{\bf Theorem \ref{reduced_charp}}
{\em Let $\Char k = p > 0$. Then the  special 
 fibre of $\Mloc$ is reduced.}

We give an outline of the proof:

We embed the special fibre $\Mlocs$ of the local model 
into the affine flag variety
 $\cF = SL_n(k\xT)/B$. Set-theoretically $\Mlocs$ is
a union of Schubert varieties. 

Further consider the special fibres $\Mlocs_0$ resp. $\Mlocs_{0,\kappa}$
of the local models of type $\{0\}$ resp. $\{0,\kappa\}$.
We can embed them in $SL_n(k\xT)/P^0$ resp. $SL_n(k\xT)/P^{0,\kappa}$,
where $P^I$, $I\subseteq\{0,\dots,n-1 \}$, is the stabilizer
of the lattice chain corresponding to $I$.
Denote the inverse images under the canonical projections
in $\cF$ by $\Mloct_0$ resp.
 $\Mloct_{0,\kappa}$.

Obviously, 
 $$ \Mlocs = \bigcap_{\kappa =1}^{n-1} 
              \Mloct_{0,\kappa}. $$

Now, $\Mloct_0$ is invariant under the action of the Iwahori
subgroup, and thus is set theoretically a union of Schubert varieties.
But we also know that it is a smooth, connected scheme, so in particular,
is reduced and irreducible. Thus it is a Schubert variety
in $\cF$. 

Furthermore, it can be shown that the $\Mlocs_{0,\kappa}$
essentially are so-called varieties of circular complexes
(see \cite{MT} resp. section \ref{equ_mlocmk}): locally, they have the form
$$ \{ (X,Y) \in \Mat_{N'}(k) \times \Mat_{N'}(k); \
      XY = YX = 0 \} \times \A^N . $$
In fact, locally all the $\Mloc_I$ can be interpreted
as spaces of certain homomorphisms (up to a product with affine space),
see theorem \ref{descr_as_homom}.
Now by
the results of Strickland \cite{Strickland} respectively of  
Mehta and Trivedi \cite{MT}, the $\Mlocs_{0,\kappa}$
are reduced.  Hence the $\Mloct_{0,\kappa}$ are reduced as well
and thus they are
unions of Schubert varieties even scheme-theoretically.

To prove the theorem, we apply the technique of Frobenius splittings.
As intersections and unions of compatibly split subvarieties are
split again, and split schemes are reduced, the theorem follows
from (cf. corollary \ref{schubvar_fsplit}):

{\bf Theorem.}
{\em \ The Schubert variety $\Mloct_0$ is Frobenius split,
and all Schubert subvarieties of $\Mloct_0$
 are simultaneously compatibly split.}

The corresponding theorem is well known for the finite dimensional 
flag variety
(we will recall this briefly in section \ref{fsplit_classic}).
Mathieu proved a similar theorem in the context of Kac-Moody algebras
(cf. \cite{M}).

Once one knows that the special fibre of the local model is reduced,
it is not very difficult to show that the local model
itself is flat over $O$.
Namely, an explicit calculation yields that the generic points of
the irreducible components of the special fibre can be lifted to
the generic fibre (see proposition \ref{Ux_smooth}).

The consideration of the special fibre of the standard
local model leads to the following question on incidence
varieties of flag varieties.
Fix $n>0 $ and a partition $\underline{r}$ of $n$. 
Let $V_0, \dots, V_\ell$ be vector spaces of dimension $n$
over some field $k$. Choose 
$(\varphi_{ij})_{ij} \in \prod_{i,j} \hom(V_i,V_j)$ and define
$$ X:= \{ (F_i)_i \in \prod_{i=0}^\ell \Flag_{\underline{r}}(V_i); \
      \varphi_{ij}(F_i) \subseteq F_j \text{ for all } i,j \}, $$
where $\Flag_{\underline{r}}(V_i)$ is the flag variety
of flags of type $\underline{r}$ in $V_i$.
What are the singularities of the scheme $X$? 

The most interesting question that remains open at the moment
is what can be
said about other groups, especially for the symplectic group.
Of course the approach of embedding the special fibre of
the local model in an affine flag variety should work as well.
For proving the reducedness, 
the main problem is to establish the analogue of the 
results of Strickland resp. Mehta and Trivedi
that we cited in the case of $GL_n$. Some results in this direction
are already available (cf. \cite{DP}, \cite{CN}),
but this does not seem sufficient to get the proof started.
Again, it is not difficult to show that
the local model is flat, 
once one knows that the special fibre is reduced.

Finally, it is a pleasure to acknowledge the help I received
from several people with this work. First of all, I am very grateful
to M. Rapoport who initiated this work and introduced me 
into this area of mathematics. His mathematical advice
as well as his encouragement and steady interest
in my work were extremely helpful to me. Furthermore, I would like
to thank O. B\"ultel, T. Haines, S. Orlik and T. Wedhorn for
many useful discussions, and T. Wedhorn again for making a lot
of valuable remarks on this text.

\chapter{Frobenius Splittings}

In this section we give the
relevant definitions and collect some basic facts about 
Frobenius splittings. We mostly follow the article \cite{MR} of
Mehta and Ramanathan;
confer also Ramanathan's article \cite{Ram}.

\section{Definition}

Let $k$ be an algebraically closed field of characteristic $p > 0$.
Let $X$ be a $k$-scheme of finite type. 
Denote by $X'$ the base change of $X$ with 
respect to the Frobenius morphism $\Spec k \lto \Spec k$.
The relative Frobenius morphism $F: X \lto X'$ gives us 
a homomorphism $\cO_{X'} \lto F_\ast \cO_X$ of $\cO_{X'}$-modules.

\begin{Def} 
i) The scheme $X$ is called Frobenius split 
(or $F$-split), if the homomorphism $\cO_{X'} \lto F_\ast \cO_X$
admits a section. Such a section is called a splitting.

ii) Let $\sigma: F_\ast \cO_X \lto \cO_{X'}$ be a splitting.
A closed subscheme $Y \subseteq X$ with sheaf of ideals $\cI$
is called compatibly $\sigma$-split (or simply compatibly
split) if 
 $\sigma(F_\ast\cI)\subseteq \cI_{Y'}$.
\end{Def}

If $Y \subseteq X$ is compatibly $\sigma$-split, then $\sigma$
induces a splitting of $Y$.

\begin{lem} \label{split_inters}
Let $\sigma: F_\ast \cO_X \lto \cO_{X'}$ be a splitting.

i) If $Y_1, Y_2 \subseteq X$ are compatibly $\sigma$-split, then
 $Y_1 \cap Y_2$ and $Y_1 \cup Y_2$ are compatibly $\sigma$-split.

ii) If $Y = Y_1 \cup \cdots \cup Y_n \subseteq X$ is the decomposition into
irreducible components and $Y$ is compatibly $\sigma$-split, then
 $Y_1, \dots, Y_n$ are compatibly $\sigma$-split. \qed
\end{lem}

The following proposition is a trivial consequence of
the definition, but it will be very important for us.

\begin{stz} \label{splitred}  If $X$ is $F$-split, then it is reduced.
\end{stz}

{\em Proof.} If $\cO_{X'} \lto F_\ast \cO_X$ has a section, then
it must be injective. \qed

\begin{stz} \label{split_dirimage}
Let $f : Z \lto X$ be a proper morphism of algebraic varieties over $k$.
Assume $f_\ast \cO_Z = \cO_X$.

i) If $Z$ is $F$-split, then $X$ is also $F$-split.

ii) If $Y \subseteq Z$ is a closed subvariety which is compatibly split,
then its image $f(Y)$ is compatibly split in $X$. \qed
\end{stz}

{\em Proof.} Let $\sigma: F_\ast \cO_Z \lto \cO_{Z'}$ be a splitting.
Since the Frobenius morphism commutes with any morphism,
we have $f_\ast F_\ast \cO_Z = F_\ast f_\ast \cO_Z = F_\ast \cO_X$,
hence $f_\ast \sigma$ is a splitting of $X$.

Now let $Y \subseteq Z$  be compatibly $\sigma$-split. 
Let $I \subseteq \cO_Z$
(resp. $J \subseteq \cO_X$) be the ideal sheaf of $Y$
(resp. $f(Y)$). Then $f_\ast I = J$ (cf. \cite{MR}, Lemma 2),
and it follows that $f_\ast\sigma (F_\ast J) = J$.
\qed

\section{A Criterion for Splitting}

Now let $X$ be a smooth projective variety of dimension $n$ over $k$.
To find a splitting of $X$, it is enough to find
a homomorphism $F_\ast \cO_X \lto \cO_{X'}$, such that the composite
$ \cO_{X'} \lto F_\ast \cO_X \lto \cO_{X'}$ is non-zero on the fibre 
at a single point (since any homomorphim $\cO_{X'} \lto \cO_{X'}$
is a constant in $k$). Such a homomorphism, which is a splitting
up to a constant,
will also be called a splitting.

So we want to understand the global sections of 
 $\underline{\hom}_{\cO_{X'}}(F_\ast \cO_X, \cO_{X'}) = (F_\ast\cO_X)^\ast$.

Denote by $\omega_X$ the canonical bundle of $X$. Serre duality gives 
a correspondence between global sections of $\omega_X^{1-p}$ and
global sections of $(F_\ast\cO_X)^\ast$. In fact, using the Cartier 
operator, one can give a natural isomorphism 
 $F_\ast \omega_X^{1-p} \lto (F_\ast\cO_X)^\ast$, which can be written down
explicitly in terms of local coordinates on $X$.
Analyzing this isomorphism, one arrives at the following criterion
for compatible splittings:

\begin{stz} {\rm ({\cite[Prop. 8]{MR}})} \label{crit_split}
Let $X$ be a smooth projective variety of dimension $n$ over $k$.
Let $Z_1, \dots, Z_n \subseteq X$ be irreducible closed subvarieties
of codimension 1 such that for any subset $I \subseteq \{1,\dots,n\}$
the scheme-theoretic intersection $Z_I = \bigcap_{i\in I} Z_i$
is reduced and irreducible of codimension $\# I$.
Let $P = \bigcap_{i=1}^n Z_i$.

Further, suppose that there exists a global section $s$ of 
 $\omega_X^{-1}$
such that $\div\ s = Z_1 + \cdots Z_n + D$, where $D$ is an effective
divisor with $P \notin {\rm supp }\ D$. 

Then the section $s^{p-1} \in H^0(X,\omega_X^{1-p})$ gives a splitting
of $X$ which compatibly splits all the $Z_i$. \qed
\end{stz}

\section{Frobenius Splitting of Classical Schubert Varieties}
\label{fsplit_classic}

In this section we want to collect some facts about
$F$-splittings for classical Schubert varieties. We will not need
them later, but it will become clear that everything we want to do for
the affine flag variety has an analogue in the classical case.

Let $G$ be a semisimple algebraic group over $k$, and choose
a maximal torus $T$ and a Borel subgroup $B$ of $G$ which contains $T$.
We have the following theorem (\cite{Ram}):

\begin{thm} {\rm (Ramanathan)}
All Schubert varieties in $G/B$ are simultaneously
compatibly $F$-split.
\end{thm}

By lemma \ref{split_inters} and proposition \ref{splitred},
this immediately gives the following corollary. 

\begin{kor}
Arbitrary intersections of unions of Schubert varieties are reduced.
\qed
\end{kor}

We give a sketch of the proof of the theorem.

Take a reduced expression $w_0 = s_{\alpha_1} \cdots s_{\alpha_r}$ of
the longest element of the Weyl group.
Let $w_i =  s_{\alpha_1} \cdots s_{\alpha_i}$, and denote by $X_i$ the
corresponding Schubert variety. Then $X_r = G/B$.

We have the Demazure varieties $Z_i$ 
which are smooth varieties of dimension $i$ 
and maps $\psi_i: Z_i \lto X_i$. Inside $Z_i$ we have
 $i$ divisors $Z_{i1}, \dots, Z_{ii}$ (cf. \cite{Ram}).
Denote the sum of these divisors by $\partial Z_i$.

Now, we want to show that $Z_r$ is Frobenius split, and that
the $Z_{rj}$ are simultaneously compatibly split, by applying
the criterion \ref{crit_split} above. We need the following lemma.

\begin{lem}{\rm (\cite[Prop. 2]{Ram})}
The canonical bundle of $Z_i$ is 
 $\cO(-\partial Z_i) \tens \psi_i^\ast \sL_\rho^{-1}$,
where $\sL_\rho$ is the equivariant line bundle
associated to the character $\rho$ (= half the sum of the positive
roots). \qed
\end{lem}

As $\sL_\rho$ is a very ample line bundle on $X_r = G/B$,
 $\psi_i^\ast \sL_\rho$ has no base point. 
(See corollary \ref{pnobasept} ii).)
We can then apply the criterion cited above. (Compare
the proof of proposition \ref{demvarsplit}.)

But as $\psi_\ast \cO_{Z_r} = \cO_{G/B}$, it follows
from proposition \ref{split_dirimage}
that $G/B$ is Frobenius split, and that the images of
the $Z_{rj}$ are compatibly split. In particular, all
Schubert subvarieties of codimension 1 are compatibly split,
and using lemma \ref{split_inters} one can see that all
Schubert varieties are compatibly split. (Compare the
proof of corollary \ref{schubvar_fsplit}.) \qed

\chapter{The Affine Flag Variety}

Denote by $k$ an algebraically closed field.

\section{Definition and Basic Properties}

In this section, we follow the article \cite{BL} of
Beauville and Laszlo quite closely. But we do not restrict
ourselves to the case of a ground field of characteristic 0.
Indeed, we are 
especially interested
in the case $\Char k = p > 0$. 
We consider $GL_n(k\xT)$ as an ind-scheme over $k$ 
in the following way:
Define 
$$ G^{(N)}(R) = \{ g(z) \in GL_n(R\xT);\ g(z) \text{ and } 
            g(z)^{-1} \text{ have poles of order } \le N \}.$$
This is an (infinite dimensional) $k$-scheme, and we have
 $GL_n = \limind G^{(N)}$
(for further details see \cite{BL}, for example). Furthermore
let $B$ denote the standard Iwahori subgroup. It is an (infinite dimensional)
scheme over $k$. The fppf quotient
 $GL_n(k\xT)/B$ is a $k$-ind-scheme.

Similarly, we have the ind-scheme $\cF := SL_n(k\xT)/B$. 
(By abuse of notation, we denote
the Iwahori subgroup of $SL_n$ by $B$ as well). $\cF$ is
called
the affine flag variety. We will describe the ind-structure 
of $\cF$ more explicitly later.

We want to identify $\cF$ with a space of lattice chains.
First, we recall some definitions.

\begin{Def} Let $R$ be a $k$-algebra. A lattice in $R\xT^n$
is a sub-$R\xt$-module $\sL$ of $R\xT^n$ which is projective of
rank $n$, and such that $\sL \tens_{R\xt} R\xT = R\xT^n$.
Equivalently, we can say that a lattice is a sub-$R\xt$-module $\sL$
of $R\xT^n$, such that $t^N R\xt^n \subseteq \sL \subseteq t^{-N} R\xt^n$
for some $N$, and such that the $R$-module $t^{-N} R\xt^n/\sL$ 
is projective. 
\end{Def}

\begin{Def} Let $R$ be a $k$-algebra.  
A sequence $\sL_0 \subset \sL_1 \subset \dots \subset
  \sL_{n-1} \subset t^{-1} \sL_0$
of lattices in $R\xT^n$ is called a complete lattice chain, if 
 $\sL_{i+1}/\sL_i$ is a locally free $R$-module of rank 1 for all
 $i$.
\end{Def}

\begin{stz}
We have a functorial isomorphism ($R$ a $k$-algebra)
\begin{diagram}
 (GL_n(k\xT)/B)(R) & \rTo^\cong & \{ \text{complete lattice chains in } 
      R\xT^n \}.
\end{diagram}
\end{stz}

{\em Proof.}
Of course, the morphism is given by $\overline{g} \mapsto g \cdot
(R\xt^n, t^{-1}R\xt \oplus R\xt^{n-1}, \dots t^{-1}R\xt^{n-1} \oplus
R\xt)$.
To prove the proposition, one shows that Zariski-locally on
$\Spec R$, every lattice chain is of the form
$g \cdot
(R\xt^n, t^{-1}R\xt \oplus R\xt^{n-1}, \dots t^{-1}R\xt^{n-1} \oplus
R\xt)$.
This is done in \cite[App. to chapter 3]{RZ}. 
\qed

\begin{Def}
Let $r \in \Z$.
A lattice $\sL \subseteq R\xT^n$ is called $r$-special,
if $\bigwedge^n \sL = t^r R\xt$ (as submodule of
 $\bigwedge^n R\xT^n = R\xT$). A (complete) lattice chain $(\sL_i)_i$ is
called $r$-special, if $\sL_0$ is $r$-special.
\end{Def} 

\begin{stz}
Fix $r \in \Z$. Then we have a functorial isomorphism
\begin{diagram}
\cF(R) & \rTo^{\cong \ \ } 
   & \{ \text{$r$-special complete lattice chains in } 
    R\xT^n \}.
\end{diagram}
\end{stz}

{\em Proof.} 
The morphism is given by
$$ \overline{g} \mapsto g \cdot (\lambda_i)_i, $$
where
\begin{eqnarray*} 
 &&  \lambda_0 = R\xt^{n-r} \oplus tR\xt^r, 
   \lambda_1 = R\xt^{n-r+1} \oplus tR\xt^{r-1}, \\
 &&\qquad  \dots, 
   \lambda_{n-1} = t^{-1}R\xt^{n-r-1} \oplus R\xt^{r+1}.
\end{eqnarray*}
We have to show that, if $(\sL_i)_i$ is a $r$-special lattice chain,
then there exists, (fppf-)locally on $\Spec R$, an element $g \in SL_n(R\xt)$,
such that $(\sL_i)_i = g \cdot (\lambda_i)_i$.

Now by the proposition above there exists, locally on $\Spec R$
(even Zariski-locally),
 $g' \in GL_n(R\xt)$, such that
 $(\sL_i)_i = g' \cdot (\lambda_i)_i$. As $\sL_0$ is $r$-special,
we have 
\begin{diagram}
 \bigwedge^n g' : t^r R\xt = \bigwedge^n \sL_0 & \rTo^\cong &
   \bigwedge^n \lambda_0 = t^r R\xt,
\end{diagram}
so $\det(g') \in R\xt^\times$. But then we can clearly find
an element $g \in SL_n(R\xt)$,
such that $(\sL_i)_i = g \cdot (\lambda_i)_i$.
\qed

\begin{kor}
If $(\sL_i)_i$ is $r$-special, then $\sL_i$ is $(r-i)$-special,
for all $i$. \qed
\end{kor}

{\bf Remark.} It is known that $\cF$ is reduced if 
$\Char k = 0$ (cf. \cite{BL}).
It is easy to see that $GL_n(k\xT)/B$ is not reduced.

Let us now describe the ind-structure of $\cF$ in terms of lattice 
chains. Identify $\cF$ with the space of $0$-special
lattice chains as above. For $g \in SL_n(R\xT)$, we have
 $t^N R\xt \subset g R\xt \subset t^{-N}R\xt$
if and only if $g$ and $g^{-1}$ have poles (with respect to $t$) of
order $\le N$.
Thus we define 
$$ \cF^{(N)} = \{ (\sL_i)_i \in \cF; \ 
         t^N R\xt \subseteq \sL_0 \subseteq t^{-N} R\xt \}. $$
This is a closed subscheme of a finite flag variety
(consisting of certain flags in $t^{-N-1}R\xt/t^N R\xt$),
and $\cF = \limind \cF^{(N)}$. Even in characteristic 0,
it is not clear whether the schemes $\cF^{(N)}$ are reduced.

If $P \supseteq B$ is a parahoric subgroup of $SL_n(k\xT)$,
then of course we have the quotient $SL_n(k\xT)/P$
and can interpret it as a space of (partial) lattice chains again.
Namely $P$ is the stabilizer of a partial lattice
chain $(\lambda_i)_{i\in I}$, for some $I \subseteq \{ 0, \dots, n-1 \}$,
and the ind-scheme $SL_n(k\xT)/P$
parametrizes partial special lattice
chains $(\sL_i)_{i \in I}$.

\begin{stz} \label{projsmooth}
Let $P \supseteq B$ be a parahoric subgroup of $SL_n(k\xT)$.
Then the canonical projection $\cF \lto SL_n(k\xT)/P$
is a smooth morphism, the fibres of which are (finite dimensional)
flag varieties.
\end{stz}

{\em Proof.} It follows from the infinitesimal lifting criterion
that the morphism is smooth. It is clear that the fibres are
flag varieties. \qed

Finally, we introduce some more notations:
Denote by $S=\{s_0, \dots, s_{n-1}\}$ the set of  simple reflections. 
For $I\subseteq S$ denote by $P_I \supseteq B$ 
the corresponding parahoric subgroup. 
(If $P^I$ denotes the stabilizer of the lattice chain 
 $(\Lambda_i)_{i\in I}$, then we have
 $P^I = P_{\{0,\dots, n-1\}-I}$.)
We also write $P_i$ instead
of $P_{\{s_i\}}$. So $P_i$ is the subgroup of $SL_n(k\xT)$
that stabilizes all the $\Lambda_j$ except for $\Lambda_i$.
Denote by $W_a$ the affine Weyl group, and by $W$ the finite
Weyl group of $SL_n$.

\section{Schubert Varieties}

For $w \in W_a$, we have the Schubert cell $BwB/B \subseteq \cF$.
It is contained in some finite dimensional part of $\cF$.
Its Zariski closure (with the reduced scheme structure) is called
the Schubert variety associated to $w$ and denoted by $X_w$. 

Set theoretically, $X_w$ is the disjoint union
$$ X_w = \bigcup_{v \le w} BvB/B, $$
where $\le$ is the Bruhat order.

\section{Demazure Varieties}

Many of the definitions and results of this section can be found
in a similar form in \cite{M}. In fact, Mathieu shows that Schubert 
varieties in the affine flag variety are Frobenius split, just as we
want to do.
But the difference is that he gets the scheme structure in a
different way, using the theory of Kac-Moody algebras. To get
a relation to the local model, we need the
--- a priori different --- scheme structure
from the 'functorial approach' chosen above.
It is possible that these two scheme structures coincide 
(in characteristic zero, this has been shown by Faltings,
cf. \cite{BL}), but I do not know how to prove this.

\subsection{Definition}

Let $w \in W_a$ and let $ \tilde{w} = s_{i_1} \cdots s_{i_\ell}$ be
a reduced expression for $w$. (We write $\tilde{w}$ instead of $w$ to
indicate that the following definitions really depend on the
choice of a reduced decomposition.)

We define 
\begin{eqnarray*}
E(\tilde{w}) & := & P_{i_1} \times^B \cdots \times^B P_{i_\ell}, \\
D(\tilde{w}) & := & P_{i_1} \times^B \cdots \times^B P_{i_\ell}/B.
\end{eqnarray*}
The variety $D(\tilde{w})$ is called the Demazure variety 
corresponding to $\tilde{w}$.

If $\tilde{u} = s_{i_1} \cdots \widehat{s_{i_k}} \cdots s_{i_\ell}$ is reduced,
we have a closed immersion
\begin{eqnarray*}
&& E(\tilde{u}) \lto E(\tilde{w}), \\
&& (g_{i_1}, \dots, g_{i_{k-1}}, g_{i_{k+1}}, \dots, g_{i_\ell})
 \mapsto (g_{i_1}, \dots, g_{i_{k-1}}, 1, g_{i_{k+1}}, \dots, g_{i_\ell}),
\end{eqnarray*}
which induces a closed immersion
\begin{equation} \label{closedimm}
D(\tilde{u}) \lto
 D(\tilde{w}). 
\end{equation}

If $\tilde{w} = \tilde{u} \tilde{v}$, i. e. 
 $\tilde{u} = s_{i_1} \cdots s_{i_k}$, 
 $\tilde{v} = s_{i_{k+1}} \cdots s_{i_\ell}$, for some $k$,
 we get a canonical projection
morphism $D(\tilde{w}) \lto D(\tilde{u})$,
which is a locally trivial fibre bundle with fibre $D(\tilde{v})$.
In particular, let $\tilde{u} = s_{i_1} \cdots s_{i_{\ell-1}}$,
 $\tilde{v} = s_{i_\ell}$. Then we get a $\PP^1$-fibration
\begin{equation} \label{p1fibr}
D(\tilde{w}) \lto D(\tilde{u}). \end{equation}
The closed immersion
 $D(\tilde{u}) \lto D(\tilde{w})$ defined above is a section
of this fibration.

\begin{kor}
The Demazure variety $D(\tilde{w})$ is smooth and proper over $k$,
and has dimension $l(w)$. \qed
\end{kor}

Multiplication gives us a morphism 
 $\Psi_{\tilde{w}} : D(\tilde{w}) \lto X_w$.

\begin{stz} \label{demvar_compatible}
The morphism $\Psi_{\tilde{w}}:D(\tilde{w}) \lto X_w$ 
is proper and birational.
If $\tilde{u} = s_{i_1} \cdots \widehat{s_{i_k}} \cdots s_{i_\ell}$ 
is reduced,
these morphisms together with the closed immersion
(\ref{closedimm}) yield a commutative diagram
\begin{diagram}
D(\tilde{u}) & \rTo & X_u \\
\dInto & & \dInto \\
D(\tilde{w}) & \rTo & X_w
\end{diagram}
\end{stz}

{\em Proof.} 
This is clear, except maybe for the birationality.
But if $\tilde{w} = s_{i_1} \cdots s_{i_\ell}$, then
 $Bs_{i_1} B \times^B \cdots \times^B B s_{i_\ell} B / B$ is an open
part of $D(\tilde{w})$, and multiplication is an isomorphism
\begin{diagram}
 B s_{i_1} B \times^B \cdots \times^B B s_{i_\ell} B / B & \rTo^\cong &
BwB/B,
\end{diagram}
as is easily seen. \qed

In the Demazure variety $D(\tilde{w})$, we have $l(w)$ divisors
 $Z_1^{\tilde{w}} ,\dots, Z_{l(w)}^{\tilde{w}} $. 
These are defined inductively on the
length of $w$, as follows.

Write $\tilde{w} = \tilde{u} s_{i_\ell}$. We have a map
 $\pi: D(\tilde{w}) \lto D(\tilde{u})$, which is a $\PP^1$-fibration,
and we also have a section $\sigma: D(\tilde{u}) \lto D(\tilde{w})$
of $\pi$.

We define
\begin{equation} \label{def_divisors}
\begin{array}{rcl}
Z_i^{\tilde{w}} & := & \pi^{-1}(Z_i^{\tilde{u}}), \quad i = 1,\dots, l(w)-1 
   \\[.2cm]
Z_{l(w)}^{\tilde{w}} & := & \sigma(D(\tilde{u})).
\end{array}
\end{equation}

We denote by $Z^{\tilde{w}}$ the sum of the divisors $Z_i^{\tilde{w}}$,
and by $P^{\tilde{w}}$ the intersection $\bigcap_i Z_i^{\tilde{w}}$.

\begin{lem} Set-theoretically, $P^{\tilde{w}}$ consists only of one point.
\end{lem}

{\em Proof.}
As $\sigma$ is a section of $\pi$, we have (set theoretically):
 $$ Z^{\tilde{w}}_i \cap D(\tilde{u}) = Z^{\tilde{u}}_i,
          \quad i= 1, \dots, l(w)-1. $$ 
So if $\tilde{u} \ne 1$, we have $P^{\tilde{w}} = P^{\tilde{u}}$.
But $P^{\tilde{w}}$ is a point if $\tilde{w}$ has length 1,
so it is always only a point.
It corresponds to the natural map $D(1) \lto D(\tilde{w})$.
\qed

\begin{lem} \label{prop_zi}
The subvarieties $Z^{\tilde{w}}_i$ are smooth of codimension 1 in
 $D(\tilde{w})$. In the tangent space $T_{D(\tilde{w}),P^{\tilde{w}}}$
we have
$$ T_{Z^{\tilde{w}}_1,P^{\tilde{w}}} \cap \cdots \cap
   T_{Z^{\tilde{w}}_{l(w)},P^{\tilde{w}}} = \{ 0 \}. $$
In particular, $P^{\tilde{w}}$ is only a point even scheme-theoretically.
\end{lem}

{\em Proof.}
We prove the lemma by induction on the length of $w$.
Write $\tilde{w} = \tilde{u} s_{i_\ell}$ as above.
As the fibration $D(\tilde{w}) \lto D(\tilde{u})$ is locally
trivial, there exists a neighbourhood $V$ of $P^{\tilde{u}}$
in $D(\tilde{u})$ with a trivialization 
 $U:= \pi^{-1}(V) = \PP^1 \times V$.

Let $\Delta := T_{\PP^1,P} \subseteq T_{D(\tilde{w}),P}$.

Since $\pi$ is locally trivial with smooth fibres, the smoothness
of $Z^{\tilde{u}}_i$ implies that $Z^{\tilde{w}}_i$ is smooth as well
($i = 1, \dots, l(w)-1$). Finally, also 
 $Z^{\tilde{w}}_{l(w)} = D(\tilde{u})$ is smooth.

We have
\begin{eqnarray*}
T_{Z^{\tilde{w}}_i,P} & = & T_{Z^{\tilde{u}}_i,P} \oplus \Delta, 
  \quad  i= 1, \dots, l(w)-1 \\
T_{Z^{\tilde{w}}_{l(w)},P} & = & T_{D(\tilde{u}),P}.
\end{eqnarray*}
Thus, for $i= 1, \dots, l(w)-1$ we have
 $$ T_{Z^{\tilde{w}}_i,P} \cap T_{Z^{\tilde{w}}_{l(w)},P} = 
 T_{Z^{\tilde{u}}_i,P}, $$
and this immediately implies the lemma.
 \qed

\begin{lem} \label{codim1var_appear}
If $\tilde{v} < \tilde{w}$, $l(v) = l(w)-1$,
then $D(\tilde{v})$ (considered as a closed subscheme of
 $D(\tilde{w})$ by the embedding defined above) is one of
the $Z_i^{\tilde{w}}$. \qed
\end{lem}

{\em Proof.} More precisely, it is easy to see by induction
on $l(w)$ that $Z_j^{\tilde{w}}$ is the
variety $P_{i_1} \times^B \cdots \times^B P_{i_\ell}$
where $P_j$ is left out. Thus all the $D(\tilde{v})$ appear
as a $Z_j^{\tilde{w}}$. (Since (in general) not all
the $s_{i_1} \cdots \widehat{s_{i_j}} \cdots s_{i_\ell}$ are
reduced, not all of the $Z_i^{\tilde{w}}$
are of the form $D(\tilde{v})$.)  \qed

\subsection{The Canonical Bundle}

As above, let $w \in W_a$ be arbitrary and choose a reduced
decomposition $\tilde{w}$.

We want to describe the canonical bundle of the Demazure
variety $D(\tilde{w})$.
To do this, we first define a certain line bundle on
the Schubert variety $X_w$.

As above, we identify the affine flag variety with the space
of $r$-special complete lattice chains. Again we denote by $(\lambda_i)_i$
the standard $r$-special lattice chain.

The Schubert variety $X_w$ consists of certain lattice chains
 $(\sL_i)_i$. We can find $N > 0$, such that all lattices
occuring here lie between $t^{-N} k\xt^n$ and $t^N k\xt^n$.

Let $n_i = \dim_k \lambda_i/t^N k\xt^n$.
We get maps 
$$ \varphi_i \colon X_w \lto \Grass(t^{-N}k\xt^n/t^N k\xt^n, n_i), \
 (\sL_i)_i \mapsto \sL_i/t^N k\xt^n. $$ 
These maps yield a closed
embedding
\begin{equation} \label{emb_in_grass}
\begin{diagram}
\varphi\colon X_w & \rInto & \prod_{i=0}^{n-1} 
       \Grass(t^{-N}k\xt^n/t^N k\xt^n, n_i).
\end{diagram}
\end{equation}

Now let $L_i$ be the very ample generator of the Picard group
of $\Grass(t^{-N}k\xt^n/t^N k\xt^n, n_i)$, and define 
$$ L_w := \varphi^\ast \bigotimes_{i=0}^{n-1} L_i. $$
This line bundle does not depend on $N$.

\begin{stz} 
The line bundle $L_w$
on $X_w$
has the following properties:

i) If $u \in W_a$, such that $u \le w$, i.e. $X_u \subseteq X_w$,
then $L_u$ is the pull back of $L_w$.

ii) $L_w$ is very ample. \qed
\end{stz}

Denote by $L_{\tilde{w}}$ the pull back of $L_w$ along
the morphism $\Psi_{\tilde{w}}: D(\tilde{w}) \lto X_w$.

\begin{kor} \label{pnobasept}

i) Let $\tilde{u} < \tilde{w}$, $l(\tilde{u}) = \tilde{w} - 1$.
 The pull back of $L_{\tilde{w}}$ along the embedding
 $\sigma: D(\tilde{u}) \lto D(\tilde{w})$ is $L_{\tilde{u}}$.

ii) The line bundle $L_{\tilde{w}}$ does not have a base point.
\end{kor}

{\em Proof.} 
i) Apply part i) of the previous proposition and proposition
 \ref{demvar_compatible}.

ii)
It is easy to see that the pull back of any very ample line
bundle under a morphism $Z \lto Z'$, $Z' \ne \{ \rm{pt} \}$,
is base point free.
%
% Let $Q \in D(\tilde{w})$, and denote the image of $Q$
%in $X_w$ by $Q'$. As $L_w$ is very ample,
%  we can find an effective divisor $D$ on $X_w$, such that
% $L_w = \cO_{X_w}(D)$ and $Q' \not\in \supp(D)$. 
%But then $Q$ is not contained in 
%the support of the pull back of $D$ to $D(\tilde{w})$ and hence
%is not a base point of $L_{\tilde{w}}$.
\qed

To establish a relation between $L_{\tilde{w}}$ and the canonical
bundle of $D(\tilde{w})$, we need the following proposition.

\begin{stz} \label{deg_along_fibres}
Write $\tilde{w} = \tilde{u} s_{i_\ell}$.
 The degree of $L_{\tilde{w}}$ along the fibres of 
 $\pi: D(\tilde{w}) \lto D(\tilde{u})$ is $1$. 
\end{stz}

{\em Proof.}
Take a point $(g_1,\dots,g_{i_{\ell-1}}) \in D(\tilde{u})$ and let
$ g = g_1 \cdots g_{i_{\ell-1}}$. The fibre over this point is isomorphic
to $P_{i_\ell}/B$ and maps to $g P_{i_\ell}/B$ in $X_w$. Thus
the degree along the fibre over $(g_1,\dots,g_{i_{\ell-1}})$
equals the degree of the pull back of $L_w$
to $g P_{i_\ell}/B$.

We have a closed immersion
 $g P_{i_\ell}/B \lto \Grass(t^{-N}k\xt^n/t^N k\xt^n, n_{i_\ell})$, 
the image 
of which is 
the projective line $\PP_g$ consisting of all subspaces of 
 $t^{-N}k\xt/t^N k\xt$ (of dimension $n_{i_\ell}$) 
lying between $g \lambda_{i_\ell+1}/t^N k\xt^n$ and
 $g \lambda_{i_\ell-1}/t^N k\xt^n$ 
(note that $g \lambda_{i_\ell+1}$ and
 $g \lambda_{i_\ell-1}$  really lie between
 $t^{-N}k\xt^n $ and $t^{N}k\xt^n $,
since $g(\lambda_i)_i \in X_w$).
But in view of the lemma below, the Picard groups
of the Grassmannian and of $\PP_g$ are isomorphic
(via restriction of line bundles), and this shows that
 $L_{i_\ell}|_{\PP_g}$ has degree $1$. 

On the other hand, for $i \ne i_\ell$,
the image of $ g P_{i_\ell}/B$ under 
$\varphi_i \colon X_w \lto \Grass(t^{-N}k\xt^n/t^N k\xt^n, n_i)$ 
is just a 
point, so $ (\varphi_i^\ast L_i)|_{\PP_g} $ has degree $0$.
\qed

\begin{lem}
Let $V$ be a $k$-vector space, $n < \dim V$, $U \subseteq W \subseteq V$
subspaces, such that $\dim W = n+1$, $\dim U = n-1$.
Denote by $\PP$ the projective line inside $\Grass(V,n)$
consisting of all subspaces of $V$ lying between $U$ and $W$.
Then restriction of line bundles gives an isomorphism
of the Picard groups of $\Grass(V,n)$ and $\PP$. \qed
\end{lem}

We cite the following lemma from \cite{Ram}.

\begin{lem}  {\rm (Ramanathan)}
Let $\pi:X \lto Y$ be a $\PP^1$-bundle with $X$ and
 $Y$ smooth varieties. Let $\sigma: Y \lto X$ be a section,
 $D$ the divisor $\sigma(Y)$ in $X$ and $L_D$ the line bundle
 $\cO_X(D)$ corresponding to the divisor $D$.

i) The relative canonical bundle $\omega_{X/Y} = 
     \omega_X \tens \pi^\ast \omega_Y^{-1}$ is isomorphic to 
 $L_D^{-2} \tens \pi^\ast \sigma^\ast L_D$.

ii) If $L$ is any line bundle on $X$ whose degree along the fibres of
 $\pi$ is 1, then 
 $\omega_{X/Y} \cong L_D^{-1} \tens 
                (L^{-1} \tens \pi^\ast\sigma^\ast L). $ \qed
\end{lem}

Now, the following description of the canonical bundle of 
 $D(\tilde{w})$ is a purely formal consequence of our definitions,
the corollary (part i)) and the proposition above.

\begin{stz} \label{canbundle}
The canonical bundle of $D(\tilde{w})$ is
 $$ \omega_{D(\tilde{w})} = \cO(-Z^{\tilde{w}}) \tens L_{\tilde{w}}^{-1}. $$
\end{stz}

{\em Proof.}
We do induction on the length of $w$.
If $l(w) = 1$, then $D(\tilde{w}) = P_i/B \cong \PP^1$.
As $\cO(-Z^{\tilde{w}}) \tens L_{\tilde{w}}^{-1}$ has degree $-2$,
it is the canonical bundle.

Now let $l(w)>1$. Again, we write $\tilde{w} = \tilde{u} s_{i_\ell}$,
 $\pi:D(\tilde{w})\lto D(\tilde{u})$, 
 $\sigma: D(\tilde{u})\lto D(\tilde{w})$ 
(see (\ref{p1fibr}), (\ref{closedimm})).
First, note that 
$$ \omega_{D(\tilde{w})} = \pi^\ast \omega_{D(\tilde{u})} 
   \tens \omega_{D(\tilde{w})/D(\tilde{u})}.$$
Let $D$ be the divisor $D(\tilde{u})$ in $D(\tilde{w})$,
and denote by $L_D$ the associated line bundle.
As $L_{\tilde{w}}$ has degree 1 along the fibres of $\pi$, 
the preceding lemma 
gives us
$$ \omega_{D(\tilde{w})/D(\tilde{u})} \cong 
   L_D^{-1} \tens (L_{\tilde{w}})^{-1}
     \tens \pi^\ast\sigma^\ast L_{\tilde{w}}. $$
By induction hypothesis, 
$$ \omega_{D(\tilde{u})} = \cO(-Z^{\tilde{u}}) \tens 
                              (L_{\tilde{u}})^{-1}.$$
But by construction (see (\ref{def_divisors})),
$$ \pi^\ast Z^{\tilde{u}} = Z^{\tilde{w}} - D. $$
Finally, we have 
$$ \sigma^\ast L_{\tilde{w}} =L_{\tilde{u}}.$$ 
Thus we get
\begin{eqnarray*}
\omega_{D(\tilde{w})} &  = &  
\pi^\ast \omega_{D(\tilde{u})} \tens \omega_{D(\tilde{w})/D(\tilde{u})} \\
& = & \pi^\ast(\cO(-Z^{\tilde{u}}) \tens (L_{\tilde{u}})^{-1})
     \tens \omega_{D(\tilde{w})/D(\tilde{u})} \\
& = & \cO(-Z^{\tilde{w}}) \tens L_D \tens \pi^\ast  
     (L_{\tilde{u}})^{-1}
     \tens  L_D^{-1} \tens (L_{\tilde{w}})^{-1}
     \tens \pi^\ast\sigma^\ast L_{\tilde{w}} \\
& = & \cO(-Z^{\tilde{w}}) \tens \pi^\ast (L_{\tilde{u}})^{-1}
     \tens (L_{\tilde{w}})^{-1}
     \tens \pi^\ast L_{\tilde{u}} \\
& = & \cO(-Z^{\tilde{w}}) \tens (L_{\tilde{w}})^{-1}.
\end{eqnarray*}
This is precisely what we wanted.
\qed

\section{Normal Schubert Varieties are $F$-split}

Now assume that our algebraically closed field $k$ has 
characteristic $p > 0$.

\begin{stz} \label{demvarsplit}
The Demazure variety $D(\tilde{w})$ admits a Frobenius
splitting which compatibly splits
all the divisors $Z_i^{\tilde{w}}$. 
\end{stz}

{\em Proof.} We want to apply the
criterion of Mehta and Ramanathan
for Frobenius splitting (proposition \ref{crit_split}).

By lemma \ref{prop_zi}, the divisors $Z_i^{\tilde{w}}$ satisfy
the necessary conditions. 
By proposition \ref{canbundle} 
the canonical bundle on $D(\tilde{w})$ 
is $\cO(-Z) \tens L_{\tilde{w}}^{-1}$. Since
the point $P$ ( = $\bigcap_i Z_i^{\tilde{w}}$) is not
a base point of $L_{\tilde{w}}$ by corollary \ref{pnobasept} ii),
we can find a global section $t$ of $L_{\tilde{w}}$, such that $P$ is not 
contained in the support of the (effective) divisor $\div(t)$. 
But then we obviously get a global section $s$ of 
 $\omega_{D(\tilde{w})}^{-1}$, such that 
$$ \div\ s = Z_1^{\tilde{w}} + \cdots Z_n^{\tilde{w}} + \div(t), $$
and this shows that $D(\tilde{w})$ is Frobenius split, and
that all the $Z_i^{\tilde{w}}$ are compatibly split.
\qed

\begin{kor} 
Assume that $X_w \subseteq \cF$ is a normal Schubert variety.
Then $X_w$ is $F$-split, and all Schubert varieties of
codimension 1 in $X_w$ are simultaneously compatibly split.
\end{kor}

{\em Proof.} Use proposition \ref{split_dirimage} and lemma 
 \ref{codim1var_appear}. \qed

\begin{lem}
Let $v, w  \in W_a$, $ v<w$, $l(v) = l(w)-2$. Then there exist
(precisely) two elements $v' \in W_a$ with $v < v' < w$. 
\end{lem} 

{\em Proof.} This is proved for a finite Weyl group
in \cite{Dixmier}, lemma 7.7.6. The same proof applies for the affine
Weyl group.\qed

\begin{kor} \label{schubvar_fsplit}
Assume that $X_w \subseteq \cF$ is a normal Schubert variety.
Then $X_w$ is $F$-split, and all Schubert subvarieties in $X_w$
are simultaneously compatibly split.
\end{kor}

{\em Proof.} This follows from the previous corollary
by induction on the codimension of the Schubert subvariety: 
The corollary says that the codimension 1 Schubert subvarieties
are compatibly split. Now assume we knew the codimension $i$
Schubert subvarieties to be compatibly split and take one
of codimension $i+1$, say $Y$. Of course, there
is a Schubert variety $X' \subseteq X$ of codimension $i-1$
which contains $Y$. By the lemma above, we find Schubert
varieties $X'_1$ and $X'_2$ in $X'$, such that $X'_1 \ne X'_2$ and
 $ Y \subset X'_1 \cap X'_2 $. Because of dimension reasons,
 $Y$ must be an irreducible component of $X'_1 \cap X'_2 $
and thus is compatibly split as well by lemma 
\ref{split_inters}. \qed

{\bf Remark.} It should be expected that all Schubert varieties
are normal. This is true in the context of Kac-Moody algebras,
cf. \cite{M}.

In any case, it is clear that the theorem holds for all Schubert varieties
that are embedded in a normal Schubert variety. In fact, it follows from
the above that all those Schubert varieties are themselves normal
(see the next section).

\section{Consequences}
\label{sing_schub_var}

The fact that normal Schubert varieties are Frobenius split
 allows one
to draw conclusions about their singularities.
First of all, we have

\begin{stz}
Let $X_w \subseteq \cF$ be a normal Schubert variety.
Then all Schubert subvarieties of $X_w$ are normal.
\end{stz} 

{\em Proof.} Let $X_u$ be a Schubert subvariety of $X_w$.
Then $X_u$ is Frobenius split, hence to prove that $X_u$ is normal
it is enough to find
a normal variety $D$ and a surjection $D \lto X_u$ with connected fibres
(see \cite{MS}). Clearly we can take the Demazure variety $D(\tilde{u})$
for some reduced expression $\tilde{u}$ of $u$. \qed

Furthermore, we show that normal Schubert varieties have
rational singularities. More precisely:

\begin{thm}
Let $X_w \subseteq \cF$ be a normal Schubert variety, and choose
a reduced expression $\tilde{w}$.
Then the morphism $\Psi_{\tilde{w}} : D(\tilde{w}) \lto X_w$ 
is a rational resolution.
In particular, $X_w$ is Cohen-Macaulay.
\end{thm}

{\em Proof.} The proof of theorem 4 in \cite{Ram} in principle works
in our situation as well. The main steps of the proof
are the following:

We have to show that 
\begin{enumerate}
\item $\Psi_\ast \cO_{D(\tilde{w})} = \cO_{X_w}$,
\item $R^q \Psi_\ast \cO_{D(\tilde{w})} = 0$ for $q>0$,
\item $R^q \Psi_\ast \omega_{D(\tilde{w})} = 0$ for $q>0$.
\end{enumerate}

A morphism that satisfies the first two conditions is called trivial.

The first point is clearly fulfilled since $X_w$ is normal
and $\Psi_{\tilde{w}}$ is proper and birational.

The third point follows from the Grauert-Riemenschneider theorem
for Frobenius split varieties (see \cite{MvK}).

The second point will be proved by induction on
the length of $\tilde{w}$. We will follow the proof of
theorem 4 in \cite{Ram}.

Write $\tilde{w} = \tilde{u} s_{i_\ell}$. The following diagram
is cartesian:
\begin{diagram}
D(\tilde{w}) & & \rTo^{\Psi_{\tilde{w}}} & & X_w \\
\dTo         & &                         & & \dTo^{\pr} \\
D(\tilde{u}) & \rTo^{\Psi_{\tilde{u}}} & X_u & \rTo^\pr & \pr(X_w)
\end{diagram}
Here $\pr$ denotes the projection $\cF \lto SL_n(k\xT)/P_{i_\ell}$.

Since the composition of trivial morphisms is trivial
and triviality is stable under flat base change, 
it is enough to show that $\pr : X_u \lto \pr(X_w) (=\pr(X_u))$
is trivial.
To show this, we apply the following criterion of Kempf (compare
\cite{Ram}, prop. 3).

\begin{lem}
Let $f: X \lto Y$ be a proper morphism
of algebraic varieties and let $L$ be an ample line bundle on $Y$ 
such that $f_\ast \cO_X = \cO_Y$ and  $H^q(X,f^\ast L^n)=0$ for $q>0$
and $n$ large. Then $f$ is trivial. \qed
\end{lem}

Consider the embedding (\ref{emb_in_grass}):
\begin{diagram}
\varphi\colon X_u & \rInto & \prod_{i=0}^{n-1} 
       \Grass(t^{-N}k\xt^n/t^N k\xt^n, n_i).
\end{diagram}
We denote the line bundle $\varphi_i^\ast L_i$ on $X_u$ by $L_{u,i}$.
Note that the morphism 
$\varphi_i : X_u \lto \Grass(t^{-N}k\xt^n/t^N k\xt^n, n_i)$ factors 
through $\pr(X_u)$. The pull-back $L'_i$ of the line bundle $L_i$
to $\pr(X_u)$
is very ample.

We want to apply the lemma above to the morphism 
 $\pr : X_u \lto \pr(X_w)$ and the line bundle $L'_i$. 
It can be shown as in the proof of theorem 2 in \cite{Ram}
that  $H^q(X_u, L_{u,i}) = 0$ for $q>0$,
and thus the hypothesis of the lemma is fulfilled.
\qed

The theorem holds also in characteristic 0. Probably this can
be derived from the results in positive characteristic by
some kind of continuity argument, but I have not thoroughly checked
this. In any case one can use the results of Mathieu in \cite{M}
since it is known that in characteristic 0 the affine flag manifolds 
coincide.

\chapter{The Local Model}

\section{Definition of the Standard Local Model}
\label{def_Mloc}

Let $O$ be a complete discrete valuation ring with perfect 
residue class field.
% of characteristic $p > 0$.
Let $\pi$ be a uniformizer of $O$ and let $k$ be an algebraic
closure of the residue class field of $O$.

Denote the quotient field of $O$ by $K$. Let $e_1, \dots, e_n$
be the canonical basis of $K^n$.

Let $\Lambda_i$, $0 \le i \le n-1$, be the free $O$-module of rank 
$n$ with basis 
 $e^i_1 := \pi^{-1} e_1, \dots,e^i_i := \pi^{-1} e_{i},e^i_{i+1} :=e_{i+1}, 
   \dots, e^i_n := e_n$.
This yields a complete lattice chain
\begin{diagram}
\cdots & \rTo & \Lambda_0 & \rTo & \Lambda_1 
& \rTo & \cdots & \rTo & \Lambda_{n-1} & \rTo & \pi^{-1}\Lambda_0 & \cdots
\end{diagram}

Fix a dominant minuscule cocharacter $\mu = (1^r,0^{n-r})$ of $GL_n$
(with respect to the torus of diagonal matrices and the Borel subgroup
of upper triangular matrices).
% (We may assume that $r \le \frac{n}{2}$.)

Furthermore choose $I = \{ i_0 < \cdots < i_{m-1} \} 
 \subseteq \{0, \dots, n-1\}$.

The standard local model $\Mloc_I$ is the $O$-scheme 
representing the following functor (cf. \cite{RZ}, definition 3.27):

For every $O$-scheme $S$, $\Mloc_I(S)$ is the set of isomorphism classes
of commutative diagrams 
\begin{diagram}
\Lambda_{i_0, S} & \rTo & \Lambda_{i_1,S} & \rTo & \cdots 
& \rTo & \Lambda_{i_{m-1},S} & \rTo^\pi & \Lambda_{i_0,S} \\
\uInto              &      & \uInto        &      &   &&\uInto & &\uInto   \\
\cF_0     & \rTo & \cF_1  & \rTo & \cdots & \rTo & \cF_{m-1}
  &\rTo & \cF_0
\end{diagram}
where $\Lambda_{i,S}$ is $\Lambda_i \tens_{O} \cO_S$, and where 
the $\cF_\kappa$ are locally free $\cO_S$-submodules of rank $r$ which
Zariski-locally on $S$ are direct summands of $\Lambda_{i_\kappa,S}$. 

It is clear that this functor is indeed representable. In fact, 
 $\Mloc_I$ is a
closed subscheme of a product of Grassmannians.

We write $\Mloc := \Mloc_{\{0, \dots, n-1\}}$.
Furthermore, we will often write $\Mloc_{i_0, \dots, i_{m-1}}$
instead of $\Mloc_{\{i_0, \dots, i_{m-1} \}}$.

As the stabilizer of the complete lattice chain is
the Iwahori subgroup, and the stabilizer of a partial
lattice chain is a parahoric subgroup, we speak also
of the local model in the Iwahori case resp.
in the parahoric case.

% More generally, if $I \subseteq \{0, \dots, n-1 \}$, we can define 
% a parahoric local model $\Mloc_I$ by considering not the complete
% lattice chain, but only those $\Lambda_i$ with $i \in I$.

Note the following obvious fact:

\begin{lem} \label{spec_fibr_agree}
Let $\Mloc$ be the local model over $O$ as above,
and let ${\Mloc}'$ be the local model over $k\xt$.
Then the special fibres of $\Mloc$ and ${\Mloc}'$ are the same.\qed
\end{lem}

\section[The Local Model and the Affine Flag Variety]
        {The Standard Local Model and the Affine Flag Variety}
\label{locmod_afffm}

We can embed the special fibre $\Mlocs$ of $\Mloc$ in the affine
flag variety
(over $k$) as follows: 

We identify $\cF$ with the space of $(n-r)$-special lattice 
chains.
Let $R$ be a $k$-algebra and write 
\begin{eqnarray*} 
 &&  \lambda_0 = R\xt^{n}, 
   \lambda_1 = t^{-1} R\xt^{1} \oplus R\xt^{n-1}, \\
 &&\qquad  \dots, 
   \lambda_{n-1} = t^{-1}R\xt^{n-1} \oplus R\xt^{1}.
\end{eqnarray*}
Let $(\cF_i)_i$ be an $R$-valued point of $\Mloc$.
Then $\cF_i$ is a subspace in 
 $\Lambda_i = R^n \cong \lambda_i / t \lambda_i$.
Let $\sL_i$ be the inverse image of $\cF_i$ under
the canonical projection $\lambda_i \lto \lambda_i / t \lambda_i$.
This gives us a complete lattice chain $(\sL_i)_i$.

\begin{lem}
The complete lattice chain defined above is $(n-r)$-special. \qed
\end{lem}

Thus we have defined a point of $\cF$ and we get a closed immersion
\begin{equation}
 \ii : \Mlocs \lto \cF.
\end{equation}
Via $\ii$, $\Mlocs$ is identified with the closed subscheme
of $\cF$ consisting of those lattice chains $(\sL_i)_i$
with $\lambda_i \supseteq \sL_i \supseteq t\lambda_i$
for all $i$.

In the same way we get closed immersions $ \Mlocs_I \lto 
 SL_n(k \xT)/P^I $
for each subset $I \subseteq \{0, \dots, n-1 \}$, where $P^I$ denotes
the (elementwise) stabilizer of the lattice chain $(\lambda_i)_{i\in I}$.
It is a parahoric subgroup of $SL_n(k \xT)$.
Of course, $\Mlocs_I$ denotes the special fibre of $\Mloc_I$.

\section{The Stratification of the Special Fibre}

Consider $\Mlocs$ as a closed subscheme of $\cF$.
It is invariant under the action of the Iwahori subgroup $B$
and thus set theoretically is a union
of Schubert varieties. Thus the decomposition into
Schubert cells gives us a stratification of
the special fibre of the local model.
 
Consider the Bruhat-Tits building of $SL_n$ over 
 $k\xT$.
Identify the vertices of the standard apartment with $\Z^n/\Z$,
such that the lattice generated by 
 $t^{-r_1}e_1,\dots,t^{-r_n}e_n$ corresponds
to $(r_1, \dots, r_n)$.
Let $\omega = (\omega_1,\dots,\omega_n)$ be the base alcove, i.e. 
 $\omega_i = (1^i, 0^{n-i})$. Denote by $\tau$ the alcove
 $((1^r, 0^{n-r}), (1^{r+1}, 0^{n-r-1}), \dots, (2^r, 1^{n-r}))$.

Furthermore, recall the following definitions from
the article \cite{KR} by Kottwitz and Rapoport.

\begin{Def}
Let $x=(x_1, \dots, x_n)$ be an alcove.

i) The number $\sum(x_i) - \sum(\omega_i) = 
     \sum_j x_i(j) - \sum_j \omega_i(j)$ is independent of $i$ 
   and is called the size of $x$.

ii) We say that $x$ is minuscule if
$$ 0 \le x_i(m) - \omega_i(m) \le 1 \quad \text{for all }
   i \in \{0,\dots,n-1\},\ m \in \{1,\dots, n\}. $$

iii) We say that $x$ is $\mu$-admissible,
if $x \le w(\mu)$ for some $w \in W$, where $W$ denotes 
the finite Weyl group.
\end{Def}

In fact, an alcove is $\mu$-admissible if and only if it is minuscule
of size $r$ (see \cite{KR}, theorem 3.5).

Now let $(\cF_i)_i \in \Mloc(k)$. 
Let  $(\sL_i)_i$ be the associated lattice chain in $\cF$, as above.
Then there is a uniquely determined element $w$ in the affine
Weyl group $W_a$ (referred to as the relative position
of $(\sL_i)_i$ and $\tau$), 
and an element $b$ in the Iwahori subgroup 
 $B$, such that $(\sL_i)_i = bw\tau$.
We say that $w\tau$ is the alcove associated to the point
 $(\cF_i)_i$.

The alcoves which occur here are just the minuscule ones
(of size $r$), in other words the $\mu$-admissible alcoves.

We get
 $$ \Mlocs = \bigcup_{x\ \mu-{\rm adm.}} S_x, $$
where $S_x = BwB/B$ is the Schubert cell associated to $x=w\tau$.

If $x=w\tau$ is a $\mu$-admissible alcove, then
we denote by $l(x)$ the length of $w \in W_a$.

\begin{lem} The stratification has the following properties:

i) The strata are just the orbits of the action of 
 $B$ 
on $\Mloc$.

ii) We have $S_x \subseteq \bar{S_y}$ if and only if $x \le y$ with respect
to the Bruhat order.

iii) The dimension of $S_x$ is $\dim S_x = l(x)$.  \qed
\end{lem}

The stratum corresponding to the alcove $\tau$ only consists of one point.
This is the worst singularity of $\Mlocs$. (While it is difficult to give 
a precise meaning to the term 'worst singularity', it should be
intuitively clear what is meant: As the singularities cannot become
better under specialization, the worst singularity
has to appear in the one-point stratum $S_\tau$,
which lies in the closure of every other stratum. In other words,
if $\Mlocs$ has a certain nice property in the point $S_\tau$
(for example reducedness), then this should hold everywhere.)

We want to associate to each $\mu$-admissible alcove $x$ an open 
subset $U_x$
of $\Mloc$ which contains the corresponding stratum. 
 
Let $x= w\tau = (x_1, \dots, x_n) $ be a $\mu$-admissible alcove. 
If $(\cF_i)_i \in S_x$,
then we have
 $$ \sL_i = b \cdot    
    \left( \begin{array}{cccc}
            t^{-x_i(1)+1} & & & \\
            & t^{-x_i(2)+1} & & \\
            & & \ddots & \\
            & & & t^{-x_i(n)+1} \\
           \end{array} \right),
$$
for some $b$ in the Iwahori subgroup.
Here we think of the matrix on the right hand side as the submodule
of $k\xT^n$ generated by the column vectors 
(with respect to the canonical basis
 $e_1,\dots,e_n$).

Again let $\lambda_i = t^{-1} k\xt^i \oplus k\xt^{n-i}$.
Denote by $e^i_1 = t^{-1}e_1, \dots, e^i_n = e_n$ the
canonical basis of $\lambda_i$.

Then $\omega$ corresponds to the lattice chain $(\lambda_i)_i$,
i.e. $(1^i, 0^{n-i})$ corresponds to $\lambda_i$. As we want to
consider $\sL_i$ as a submodule of $\lambda_i$
and as $\omega_i \le x_i \le \omega_i + (1, \dots, 1)$, we indeed have
to take $-x_i(\cdot)+1$ instead of $-x_i(\cdot)$ 
as the exponent in the matrix above. 

% Since we consider all these lattices
% only up to homothety, this does not really matter, of course.

The above description of $\sL_i$ shows
that the quotient $\lambda_i/\sL_i$ is generated by those
 $e^i_j$ with $\omega_i(j) - (x_i(j) -1 ) = 1$, i.e. with
 $\omega_i(j) = x_i(j)$.

Thus the open subset of $\Mlocs$, where for all $i$ the quotient
 $\lambda_i/\sL_i$ is generated by those $e^i_j$ with
 $\omega_i(j) = x_i(j)$, contains the stratum $S_x$.

We want to define an open subset of $\Mloc$ which contains
the stratum $S_x$. Thus we consider more generally the
quotients  $\Lambda_i/\cF_i$ and define (compare
lemma \ref{spec_fibr_agree}):

\begin{Def}
Let $x = (x_1,\dots, x_n)$ be a minuscule alcove of size $r$.
Then let $U_x$ be the open subset of $\Mloc$ which consists of
all points $(\cF_i)_i$, such that for all $i$ the quotient
 $\Lambda_i/\cF_i$ is generated by those $e^i_j$ with
 $\omega_i(j) = x_i(j)$.
\end{Def}

We have

\begin{stz} \label{propsUx}
i) The stratum $S_x$ is contained in $U_x$.

ii) The open subset $U_\tau$ intersects every
stratum.

iii) The irreducible components of the special fibre $\Mlocs$ 
   are the closures
   of the $U_x \cap \Mlocs$, where $x$ is an extreme alcove, i. e.
   $x=w(\mu)$ for some $w\in W$.
\end{stz}

{\em Proof.} Part i) follows from the discussion above, and
ii) and iii) follow from the lemma. \qed

Of course, we can just as well characterize the irreducible components 
of $\Mlocs$ as the closures of the strata 
 $S_{x}$, $x\in W(\mu)$. Thus $\Mlocs$ has
 $\# W/W_\mu = \left( n \atop r \right) $ irreducible components,
where $W_\mu \subseteq W$ denotes the stabilizer of $\mu$.
A similar description can be given in the parahoric case.

Finally, we note the following lemma:

\begin{lem} \label{adm_is_adm}
Let $x = (x_1,\dots, x_n)$ be a minuscule alcove of size $r$.
Let $t_i = x_i - \omega_i$, $i = 0, \dots, n-1$.
Then 

i) For all $j$, $(t_1(j), \dots, t_n(j))$ is a cyclic
 permutation of $(1^\kappa, 0^{n-\kappa})$ for some $\kappa$.

ii) If $t_i(j)=1$ and $t_{i+1}(j)=0$, then 
 $\phi_i(e^i_j) = \pi e^{i+1}_j$, where $\phi_i$ is the map
 $\Lambda_i \lto \Lambda_{i+1}$. \qed
\end{lem}

Of course, we can do similar things in the parahoric case. 
We then have to work with 'partial alcoves' $x = (x_i)_{i\in I}$
(cf. also \cite{KR}, \S 9).
The lemma above holds then in an analogous form.

\section{The Equations of the Standard Local Model $\Mloc$}

We want to compute the equations describing the standard
local model.

\subsection{More General Schemes of Compatible Subspaces}

To establish the relation between the local model and 
certain spaces of homomorphism, as will be done in the
next section, it is useful to introduce more general
'schemes of compatible subspaces'. 

Let $m, n > 0$, $0 < r  <m$, and consider free $O$-modules
 $ \Lambda_i $, $ i = 0, \dots, m-1$, of rank $n$ with bases 
 $\underline{e}^i = (e^i_1, \dots, e^i_n)$

Take $O$-linear
maps $\phi_i : \Lambda_i \lto \Lambda_{i+1}$, $i = 0, \dots, m-1$
($\Lambda_m := \Lambda_0$).

Then we denote by $M(m,n,r,(\phi_i)_i)$ the functor
which associates to an $O$-scheme $S$
the set of isomorphism classes
of commutative diagrams 
\begin{diagram}
\Lambda_{0, S} & \rTo^{\phi_0} & \Lambda_{1,S} & \rTo^{\phi_1} & \cdots 
& \rTo & \Lambda_{m-1,S} & \rTo^{\phi_{m-1}} & \Lambda_{0,S} \\
\uInto              &      & \uInto        &      &   &&\uInto & &\uInto \\
\cF_0     & \rTo & \cF_1  & \rTo & \cdots & \rTo & \cF_{m-1} &\rTo & \cF_0
\end{diagram}
where as before $\Lambda_{i,S}$ is $\Lambda_i \tens_{O} \cO_S$, and where 
the $\cF_i$ are locally free $\cO_S$-submodules of rank $r$ which
Zariski-locally on $S$ are direct summands of $\Lambda_{i,S}$. 

Again, this functor is representable by a closed subscheme of
a product of Grassmannians.
Obviously the local models occur as a special case of this definition.

It seems to be very difficult to describe these schemes
in general. The only case we will need below is
the following:

\begin{Def}
The scheme $M(m,n,r,(\phi_i)_i)$ is called a generalized
local model, if the $\phi_i$ are diagonal matrices 
with only 1's and $\pi$'s on the diagonal
(with respect to the fixed 
bases $\underline{e}^i$) and their composition is $\pi$.
\end{Def}

The only new possibility here is that 
 $\phi_i = \id$ is allowed. Strictly speaking, 
this does not happen for the local models. On the other hand,
 as the steps where 
$\phi_i$ is an isomorphism can be neglected, these schemes still
are isomorphic to certain (parahoric) local models. Thus, in the cases
we will consider,  the new definition really is only another
notation.

Now we want to define certain open subsets 
in the generalized local model 
 $M(m,n,r,(\phi_i)_i)$. For $i = 0, \dots, m-1$, 
choose $t_i \in \{0,1\}^n$ such that 
$\# \{ j; \ t_i(j)=1 \} = r$, and such that the following
two conditions are satisfied:
\begin{bed} \label{cond}
i) For all $j$, $(t_1(j), \dots, t_n(j))$ is a cyclic
 permutation of $(1^\kappa, 0^{n-\kappa})$ for some $\kappa$.

ii) If $t_i(j)=1$ and $t_{i+1}(j)=0$, then 
 $\phi_i(e^i_j) = \pi e^{i+1}_j$.
\end{bed}

Let $U = U((t_i)_i)$ be the open subset of $M(m,n,r,(\phi_i)_i)$,
where $\Lambda_{i,S}/\cF_i$ is generated by the $e^i_j$ with
 $t_i(j) = 0$.

The open subsets $U_x$ associated to an admissible alcove 
defined above are a special case of
this definition (cf. lemma \ref{adm_is_adm}). 

% Note that for a generalized local model with some
% $\phi_i = \id$  we get more open subsets of the form
% $U((s_i)_i)$ than for the isomorphic local model where
% these steps are left out.

Finally, we state the following lemma.

\begin{lem} \label{inverse}
We have
$$ M(m,n,r,(\phi_i)_i) \cong M(m,n,n-r,(\phi_i)_i). $$ 
\end{lem}

{\em Proof.} Replace $\cF_i$ with the dual of $\Lambda_i/\cF_i$.
\qed

\subsection{Interpretation in Terms of Homomorphisms}

We will see that the open subsets $U((t_i)_i$ of a
generalized local model can
be related to certain spaces of homomorphisms between
free modules. 
This will enable us to read off the equations
of the local models almost immediately. 

As the combinatorics involved here is quite complicated,
it is probably more enlightening to look at the statement
of theorem \ref{descr_as_homom} and then to try to figure out the equations
of the parahoric local model
 $\Mloc_{\mu, \kappa}$ directly rather than to go through
the proof of the theorem. The result for this special case is stated
in section \ref{equ_mlocmk}.

Consider a generalized local model $M =M(m,n,r,(\phi_i)_i)$.

% associated to 
%  $I \subseteq \{0, \dots, n-1 \}$.
% Let $m = \# I$ and
% write $I = \{ \iota_0 < \cdots < \iota_{m-1} \}$. Denote by $\phi_i$
% the map $\diag(1^{\iota_i}, \pi^{\iota_{i+1}-\iota_i}, 1^{n-\iota_{i+1}})
%  \colon \Lambda_{\iota_i} \lto \Lambda_{\iota_{i+1}}$,
% ($\Lambda_m := \Lambda_0$, $\iota_{m} := n$)
% such that $\Mloc_I = M(m,n,r,(\phi_i)_i)$.

We will represent all the subspaces $\cF_i$ by giving $r$ generating
vectors (with respect to the basis $e^i_1,\dots,e^i_n$ of $\Lambda_i$)
which we will arrange as column vectors in a matrix.

We want to study an open subset of the form $U:=U((t_i)_i)$,
where $(t_i)_i$ satisfies the condition \ref{cond}.
Choose permutations 
 $\sigma_i \in {\cal S}_n$ such that $t_i = \sigma_i (1^r, 0^{n-r})$. 

The conditions defining $U$ can be stated in the following way:
 $U$ is the open subset of $M$
consisting of those $(\cF_i)_i$, such that 
 $\cF_i$ can be described by a matrix
$$ M_i := (b^i_{jk})_{j=1, \dots,n, \ k= 1, \dots, r} 
       :=\sigma_i \
\left( \begin{array}{cccc} 
   1     &          &        &          \\
         &     1    &        &          \\
         &          & \ddots &          \\
         &          &        & 1        \\
a^i_{11} & a^i_{12} & \cdots & a^i_{1r} \\
 \vdots  & \vdots   &        &  \vdots  \\
a^i_{n-r,1} & a^i_{n-r,2} & \cdots & a^i_{n-r,r}
\end{array} \right), $$
i. e. we have a unit matrix in the rows $\sigma_i(1), \dots, \sigma_i(r)$.
Note that the $a^i_{jk}$ are uniquely determined by $\cF_i$.

The condition that $\cF_i$ is mapped under $\phi_i$ into $\cF_{i+1}$
can be expressed in terms of matrices as follows: If $\cF_i$ 
(resp. $\cF_{i+1}$) is described by $M_i$ (resp. $M_{i+1}$),
then we must have 
$$ \phi_i M_i = M_{i+1} N_i, $$
where $N_i$ is a certain $r \times r$-matrix. In fact, $N_i$  is 
uniquely determined by $M_i$ and $M_{i+1}$ since 
certain rows of $M_{i+1}$ form a unit matrix.

% In other words, we have $U = U((\sigma_i(1^r,0^{n-r}))_i)$.
% In particular, the open subsets $U_x\subseteq \Mloc$ that appeared above 
% are of this form: 
% The defining conditions just mean
% that certain $r \times r$-minors of the matrices describing
% the $\cF_i$ are invertible.

{\bf The first lemma}

Let 
$$ S:= \{ \iota \in \{1, \dots, n\}; 
        \ t_i(\iota)=0 \text{ for all } i \}, \ s:= \# S. $$

We will show that $U$ is isomorphic to the
product of an open subset of some generalized local model
$M(m,n-s,r,(\psi_i)_i)$
with $\A^{rs}$.

For all $i$ let $\overline{\Lambda}_i := \Lambda_i / 
 \langle e^i_j; \ j \in S \rangle $. Denote the map
 $\overline{\Lambda}_i \lto \overline{\Lambda}_{i+1}$
induced by $\phi_i$ by $\psi_i$.

Let $u_i :=$ '$ t_i$ with the entries $t_i(j)$, $j\in S$,
left out', such that $u_i$ is a permutation of $(1^r, 0^{n-r-s})$. 
The conditions \ref{cond} are satisfied again (with respect to
the $\psi_i$.)

\begin{lem}
We have
$$ U \cong V \times \A^{rs}_O, $$
where
$$ V:= U((u_i)_i) \subseteq M(m,n-s,r, (\psi_i)_i). $$
\end{lem}

{\em Proof.} 
We construct a surjection $U \lto V$. 
To define the map $U \lto V$, take $(\cF_i)_i$ in $U(R)$
for some $O$-algebra $R$, and consider the associated matrices
 $M_i$.

We then define a point $(\cF_i')_i \in V$ 
in the following way:
The matrix 
associated to $\cF_i'$ (with respect to the basis
 $e^i_j$, $ j \not\in S$) is the matrix $M_i'$
obtained from $M_i$ by deleting the rows with
index in $S$. Then $M_i'$ describes a subspace
of $\overline{\Lambda}_{i,R}$.
Note that $\cF_i'$ again has rank $r$.
In this way we get a morphism $U \lto V$.

What are the fibres of this map? We have to check how 
to construct matrices $M_i$, given matrices $M_i'$.
To do this, we have to fill in the entries $b^i_{jk}$
in the $j$-th row of $M_i$
for $j \in S$. So we fix $j\in S$ and consider the 
equations which arise for these $b^i_{jk}$'s.

We distinguish the following two cases
that yield two basically different types of equations:

{\em First case: } $\phi_i(e^i_j) = \pi e^{i+1}_j.$

{\em Second case: } $\phi_i(e^i_j) = e^{i+1}_j.$

The first case occurs precisely once, say for $i = i_0$, and 
the condition that $\cF_{i_0}$ is mapped into $\cF_{i_0+1}$
gives equations
of the form ($\kappa= 1,\dots,r$)
$$ \pi b^{i_0}_{j\kappa} = 
    \text{ something depending on } 
    b^{i_0+1}_{j1}, \dots, b^{i_0+1}_{jr} \text{ and }
    b^{i_0}_{\ell, \kappa} \text{ with } \ell \not\in S. $$

The second case yields equations of the form
($\kappa= 1,\dots,r$)
$$ b^{i}_{j\kappa} = 
    \text{ something depending on } 
    b^{i+1}_{j1}, \dots, b^{i+1}_{jr} \text{ and }
    b^{i}_{\ell, \kappa} \text{ with } \ell \not\in S. $$

Thus, if we choose $b^{i_0}_{j1}, \dots, b^{i_0}_{jr} \in R$,
all $b^i_{jk}$ are uniquely determined by the equations
derived from the second case. The only question is, if
the equation given in the first case is satisfied.
But it follows from the fact that all the compositions
of the $\phi_i$ are $\pi$ that this equation is automatically
satisfied. Hence we can indeed choose the
 $b^{i_0}_{j1}, \dots, b^{i_0}_{jr}$ arbitrarily.
So 'leaving out' the $j$-th row gives us a fibration
with fibres isomorphic to $\A^r$, and as $S$ has $s$
elements,  
we see that all fibres of the map $U \lto V$ (over $R$-valued points)
are isomorphic to $\A^{rs}$.
Hence we get  $U \cong V \times \A^{rs}$.
\qed

{\bf The second lemma}

Let $U=U((t_i)_i) \subseteq M(m,n,r,(\phi_i)_i)$ be as above.

Let
$$ T:=  \{ \iota \in \{1, \dots, n\}; 
        \ t_i(\iota)=1 \text{ for all } i \}, \ t:= \# T. $$

In this second step
we show that $U = U((t_i)_i)$ is isomorphic to
a product of an open subset of some $M(m,n-t,r-t,(\xi_i)_i)$
with $\A^{t(n-r)}$. 

Let $\overline{\Lambda}_i := \Lambda_i/ \langle e^i_j; \ j\in T \rangle$.
Denote the map induced by $\phi_i$ by $\xi_i$.
Let $v_i :=$ '$ t_i$ with the entries $t_i(j)$, $j\in T$,
left out'. Thus $v_i$ is a permutation of $(1^{r-t}, 0^{n-r})$,
and the conditions \ref{cond} are satisfied again (with respect 
to the $\xi_i$).

\begin{lem}
We have 
$$ U \cong W \times \A^{t(n-r)}, $$
where
$$ W:= U((v_i)_i) \subseteq M(m,n-t,r-t, (\xi_i)_i). $$
\end{lem}

{\em Proof.}
We follow the same strategy as before: first, we define a map
 $U \lto W$, then we examine its fibres.

The map $U \lto W$ is defined as follows. Take an $R$-valued
point $(\cF_i)_i$ in $U$. The $\cF_i$ correspond
to matrices $M_i$. Consider $j\in T$. In each
 $M_i$, the $j$-th row consists of $r-1$ 0's and one 1, say in
column $k^i_j$.
Construct matrices $M_i'$ as follows: For all $j\in T$,
delete the $j$-th row from
 $M_i$, and also delete the $k^i_j$-th column. 
We get a $(n-t) \times (r-t)$-matrix, and if we denote by
 $\cF_i'$ the subspace of $\overline{\Lambda}_i$ defined
by $M_i'$, we have $(\cF_i')_i \in W$.

Now suppose we are given matrices $M_i'$ associated to
 $(\cF_i')_i \in W(R)$ and we want to define matrices $M_i$
that give an element in $U(R)$. Then for each $j\in T$ we have to choose
entries $a^i_{1,{k^i_j}},\dots,a^i_{n-r,{k^i_j}}$ such that the corresponding
subspaces are mapped into one another.
Note that $\{a^i_{1,{k^i_j}},\dots,a^i_{n-r,{k^i_j}}\}
 = \{ b^i_{\iota, k^i_j}; \ t_i(\iota) = 0 \}$.

We fix $j\in T$. Again, two different types of equations appear:

{\em First case: } $\phi_i(e^i_j) = \pi e^{i+1}_j.$

{\em Second case: } $\phi_i(e^i_j) = e^{i+1}_j.$

As before, the first case occurs precisely once, say for $i = i_0$. 
We then get 
equations of the form
$$ b^{i_0}_{\iota, k^{i_0}_j}
= \pi b^{{i_0}+1}_{\iota, k^{i_0+1}_j } 
  + \text{ something depending on } b^{i_0}_{*, k^{i_0}_j}
 \text{ and } b^{{i_0}+1}_{\iota,*},
\quad \iota=1,\dots,n.
$$
In the second case, the equations have the form
$$ b^{i}_{\iota, k^{i}_j}
= b^{{i}+1}_{\iota, k^{i+1}_j } 
  + \text{ something depending on } b^{i}_{*, k^{i}_j}
 \text{ and } b^{{i}+1}_{\iota,*},
\quad \iota=1,\dots,n.
$$
One sees easily that therefore
for $i \ne i_0$, $b^{{i}+1}_{\iota, k^{i+1}_j }$ is
determined by $b^{i}_{\iota, k^{i}_j}$ and some other
 $b$'s which are 'not related to $T$'.
This yields an isomorphism between the fibre and $\A^{t(n-r)}$:
choose the $b^{{i_0}+1}_{\iota, k^{i_0+1}_j }$ with
 $ \iota$ such that $t_{i_0+1}(\iota)=0$ arbitrarily,
then determine the other $b^{i}_{\iota, k^{i}_j}$
($\iota$ such that $t_{i_0}(\iota)=0$)
by the equations from the second case. The equations 
of the second case where $t_{i_0}(\iota)=1$
and the equations of the first case are automatically satisfied.

This proves the lemma.
\qed

{\bf The third lemma}

Again, let $U=U((t_i)_i) \subseteq M(m,n,r,(\phi_i)_i)$ be as above.
Furthermore we now assume that 
$$  \{ \iota \in \{1, \dots, n\}; 
        \ t_i(\iota)=1 \text{ for all } i \}  = \emptyset. $$
By the second lemma we can write $U$ as a product of an affine space
and a new open subset of a generalized local model 
which has this property.

We will show that in this case $U$ can be described as a space
of $m$-tuples of homomorphisms the compositions
of which are $\pi$. 

If $(\cF_i)_i \in U$, then by the definition of $U$,
for all $i$ the vectors $e^i_j$ with $t_i(j)=0$ give
a basis of $\Lambda_i/\cF_i$. Here by an abuse of notation
we write $\Lambda_i/\cF_i$ instead of $\Lambda_{i,S}/\cF_i$
for variable $S$.
With respect to these bases
 the map $\Lambda_i/\cF_i \lto \Lambda_{i+1}/\cF_{i+1}$
induced by $\phi_i$
is described by a $(n-r)\times(n-r)$-matrix $X_i$.
The columns of the matrix correspond to the $e^i_j$ with 
 $t_i(j)=0$. We denote these $j's$ by $j^i_1, \dots, j^i_{n-r}$.

 Thus we get a map
$$ \Phi: U \lto \prod_{i=0}^{m-1} \Mat_{n-r}(O). $$
We want to determine the image of this map. 

It is clear
that each $X_i$ underlies the following restriction: If
 $t_i(j)=0$ and $t_{i+1}(j)=0$, then $e^i_j$ and $e^{i+1}_j$
are part of the chosen bases of 
 $\Lambda_i/\cF_i$ and $\Lambda_{i+1}/\cF_{i+1}$. Now 
 $e^i_j$ is mapped to $e^{i+1}_j$ (or $\pi e^{i+1}_j$)
under the map $\Lambda_i \lto \Lambda_{i+1}$,
and the same must be true for 
$X_i : \Lambda_i/\cF_i \lto \Lambda_{i+1}/\cF_{i+1}$.
This observation leads to the following definition:

Let $\varepsilon^i_j= 0$, if $\phi_i(e^i_j) = e^{i+1}_j$ and
  $\varepsilon^i_j= 1$, if $\phi_i(e^i_j) = \pi e^{i+1}_j$.

Then
\begin{eqnarray*}
 && {\mathfrak X}_i := \{ X=(x_{\mu,\nu})_{\mu,\nu} \in \ \Mat_{n-r}(O),\  
        \text{such that for all } \iota = 1, \dots, n-r \\
 && \qquad          \text{with } t_{i+1}(j^i_\iota)=0, 
\text{ so } j^i_\iota = j^{i+1}_{\iota'} \text{ for some } \iota',
           \text{ we have (write $j := j^i_\iota$):} \\
 && \qquad\quad  
       x_{\mu,\iota} = \pi^{\varepsilon^i_{j}} \delta_{\mu,{\iota'}}, 
       \quad \mu = 1, \dots, n-r. \}
\end{eqnarray*}

So ${\mathfrak X}_i \subseteq \Mat_{n-r}$ is an affine subspace and
the image of $\Phi$ lies in $\prod {\mathfrak X}_i$.
 
The image of $\Phi$ is determined in the next lemma.

\begin{lem}
The map defined above gives an isomorphism
$$ U \cong \{ (X_i)_i \in \prod {\mathfrak X}_i;\ 
    X_{m-1} X_{m-2} \cdots X_0 = X_{m-2} \cdots X_0 X_{m-1}=
   \cdots = \pi \}. $$
\end{lem}

{\em Proof.}
We define the inverse map. Given a tuple of matrices $(X_i)_i$
(over some $O$-algebra $R$)
in the right hand side set, we want to define an $R$-valued point
 $(\cF_i)_i$ in $U$.

First we define maps $\alpha_i: \Lambda_{i,R} \lto R^{n-r}$.
Afterwards we want to define $\cF_i := \ker \alpha_i$.

Denote the canonical basis of $R^{n-r}$ by $f_1, \dots, f_{n-r}$.

{\bf Definition of $\alpha_i$.} If $t_i(j)=0$, say $j=j^i_\iota$,
 we must have $\alpha_i(e^i_j) = f_\iota$.

If $t_i(j)=1$, then choose $i'$ such that $t_{i'}(j)=0$. Then
 $j = j^{i'}_\iota$ for some $\iota$.
Furthermore, $\phi_{i-1} \circ \cdots \circ 
 \phi_{i'+1} \circ \phi_{i'}(e^{i'}_j)  = 
 e^i_j$ by condition \ref{cond}. Now define
$$ \alpha_i(e^i_j) := X_{i-1} \cdots X_{i'} f_\iota. $$
This is independent of the choice of $i'$ (use condition
 \ref{cond} i)). It is clear that the $\alpha_i$ are surjective.

The condition that the products of the $X_i$ are $\pi$ 
together with condition \ref{cond} ii) ensures
that in this way we  get
a commutative diagram
\begin{diagram}
\Lambda_{0, R} & \rTo & \Lambda_{1,R} & \rTo & \cdots 
& \rTo & \Lambda_{m-1,R} & \rTo & \Lambda_{0,R} \\
\dOnto^{\alpha_0}  &   & \dOnto^{\alpha_1}  &  &   &&
 \dOnto^{\alpha_{m-1}} & &\dOnto^{\alpha_0}   \\
R^{n-r}     & \rTo^{X_0} & R^{n-r}  & \rTo^{X_1} & \cdots & 
 \rTo^{X_{m-2}} & R^{n-r} & \rTo^{X_{m-1}} & R^{n-r}
\end{diagram}

Thus defining $\cF_i:=\ker \alpha_i$ we indeed get a point of $U$.
It is clear that the two maps are inverse to one another.
\qed

The three lemmas together give the following theorem:

\begin{thm} \label{descr_as_homom}
Let $U= U((t_i)_i) \subseteq M(m,n,r, (\phi_i)_i)$, and
define $s$ and $t$ as above. Then 
$$ U \cong W \times \A^{rs+t(n-r-s)}_O, $$
where 
$$ W = \{ (X_i)_i \in \prod_{i = 0}^{m-1} {\mathfrak X}_i ; \ 
   X_{m-1} X_{m-2} \cdots X_0 = X_{m-2} \cdots X_0 X_{m-1}=
   \cdots = \pi \},
$$
and where the ${\mathfrak X}_i$ are certain affine subspaces
of the space of $(n-r-s) \times (n-r-s)$-matrices
 $\Mat_{n-r-s}(O)$. \qed
\end{thm}

\subsection{The Equations of $U_\tau$}
\label{genequ}

Now we want to apply the previous theorem to find the equations
for the open subset $U_\tau \subseteq \Mloc$, which is an open 
neighborhood of the worst singularity of the local model. 
So we take $\underline{t} = ((t_i)_i) = ((1^r,0^{n-r}), 
 (1^{r+1},0^{n-r-1}), \dots, (2^{r-1},1^{n-r+1}))$.

Denote by  $A_i$  the following $(n-r)\times (n-r)$-matrix:
$$ \left(\begin{array}{cccc}
a^i_1      & 1 & & \\
a^i_2      & 0 & \ddots & \\
\vdots     &   & \ddots & 1 \\
a^i_{n-r}  &   &        & 0
\end{array}\right), $$
where the $a^i_j$ are indeterminates.

\begin{stz} \label{equ_Utau}
The open subset 
 $U_\tau$ of $\Mloc$ is isomorphic to 
 $$ \Spec O[a^i_\kappa; i=0,\dots, n-1, \kappa=1,\dots, n-r] / I, $$
where $I$ is the ideal generated by the entries of the matrices
 $$ A_{n-1} A_{n-2} \cdots A_0 -\pi,\ A_{n-2} \cdots A_0 A_{n-1} -\pi,
    \ \dots,\ A_0 A_{n-1} \cdots A_1 - \pi. $$
\end{stz}

{\em Proof.} Use theorem \ref{descr_as_homom}. We have $s=t=0$,
thus only lemma 3 is needed.  \qed

The proposition gives a handy form to write down the equations,
but in fact these are much more equations than needed.
For example in the case $n=4$, $\mu = (1,1,0,0)$,
one can describe $U_\tau$ by 6 equations, whereas
the description we have given here consists of
16 equations.

{\bf Remark.} In particular, we re-discover here the well-known
equations of the local model in the Drinfeld case (i.e. $r=n-1$). Then
the $A_i$ are $1 \times 1$-matrices, so they are just indeterminates
and we get only one equation:
$$ A_{n-1} A_{n-2} \cdots A_0 = \pi.$$

{\bf Remark.} 
If one wants to analyse the properties of $\Mloc$ starting with
the equations, it might be disturbing that the $A_i$
are not generic $(n-r) \times (n-r)$-matrices, but that some
entries are 0 resp. 1. But in fact, there is a relation to the
scheme defined in an analogous way considering generic matrices.
To make this precise, let $B_i = (b^i_{jk})_{jk}$, $i = 0, \dots, n-1$ 
be $(n-r) \times(n-r)$-matrices
of indeterminates and consider the scheme
$$ M' = \Spec O[b^i_{jk}]/B_{n-1} \cdots B_0 = B_{n-2}\cdots B_0 B_{n-1} =
    \cdots = \pi.$$

If we knew that the special fibre $\overline{M}'$ of $M'$ 
has the same dimension 
as the generic fibre, namely $(n-r)^2(n-1)$, and that
$\overline{M}'$ is Cohen-Macaulay, we could conclude
that $U_\tau \cap \Mlocs$ (and thus $\Mlocs$) 
is Cohen-Macaulay as well.

The reason is that $U_\tau \cap \Mlocs$ is a closed
subscheme of $\overline{M}'$
defined by $(n-r-1)(n-r)n$ equations,
and that the number of equations is just the difference
between the dimensions of $\overline{M}'$ and $U_\tau \cap \Mlocs$.
 
See also the remarks at the end of Faltings' article \cite{F}.

\subsection{The Open Subsets $U_x$ for Extreme Alcoves $x$}

\begin{stz} \label{Ux_smooth}
Let $x$ be an extreme alcove, i. e. $x \in W(\mu)$.
Then $U_x \cong \A_O^{r(n-r)}$.
\end{stz}

In particular, $U_x$ is a smooth open subset of $\Mloc$. It seems
reasonable to expect that the special fibre of $U_x$ 
coincides with $S_x$ for extreme alcoves $x$.

% sollte daraus folgen, dass S_x affiner Raum!!!

{\em Proof.}
For an extreme alcove $x$, we have $U_x = U((t_i)_i)$
with $t_1 = \cdots = t_n$. Thus, with notation as above, $s=n-r$, $t=r$
and the proposition follows immediately from theorem \ref{descr_as_homom}. 
\qed

\subsection{The Equations of $\Mloc_{\mu,\kappa}$}
\label{equ_mlocmk}

In this section we want to determine the equations of an
open neighbourhood of the 'worst singularity' of
 $\Mloc_{\mu,\kappa}$, $\mu, \kappa \in \{0,\dots,n-1\}$, 
 $\mu \ne \kappa$.
Recall that this is the parahoric local model, where
not the complete lattice chain, but
only the the lattices $\Lambda_\mu, \Lambda_\kappa$
are involved.
As $\Mloc_\kappa$ is simply a Grassmannian, this is the
first non-trivial case.

Obviously, we may assume that $\mu=0$
and $\kappa, r \le \frac{n}{2}$.

Let 
\begin{eqnarray*}
 U & = & U((1^r, 0^{n-r}),(0^\kappa, 1^r, 0^{n-r-\kappa}))  \\
  & = &  \{ (\cF_0, \cF_\kappa) \in \Mloc_{0,\kappa} ; \\ 
&& \cF_0 \hat{=} \left( \begin{array}{cccc} 
   1     &          &        &          \\
         &     1    &        &          \\
         &          & \ddots &          \\
         &          &        & 1        \\
a^0_{11} & a^0_{12} & \cdots & a^0_{1r} \\
 \vdots  & \vdots   &        &  \vdots  \\
a^0_{n-r,1} & a^0_{n-r,2} & \cdots & a^0_{n-r,r}
\end{array}  \right), \\
&& \cF_\kappa \hat{=} \left( \begin{array}{cccc} 
 a^\kappa_{n-r-\kappa+1,1} & a^\kappa_{n-r-\kappa+1,2} & \cdots & a^\kappa_{n-r-\kappa+1,r} \\
 \vdots & \vdots & & \vdots \\
 a^\kappa_{n-r,1} & a^\kappa_{n-r,2} & \cdots & a^\kappa_{n-r,r} \\
   1     &          &        &          \\
         &     1    &        &          \\
         &          & \ddots &          \\
         &          &        & 1        \\
a^\kappa_{11} & a^\kappa_{12} & \cdots & a^\kappa_{1r} \\
 \vdots  & \vdots   &        &  \vdots  \\
a^\kappa_{n-r-\kappa,1} & a^\kappa_{n-r-\kappa,2} & \cdots & a^\kappa_{n-r-\kappa,r}
\end{array}  \right) \}.
\end{eqnarray*}

This is an open subset of $\Mloc_{0,\kappa}$, which contains
the 'worst singularity'. 
Theorem \ref{descr_as_homom} gives the following description
of $U$:

{\em First case: $\kappa\le r$}

Let $A = (a^0_{i,j})_{i, j = 1, \dots \kappa}$, 
 $ B = (a^\kappa_{i,j})_{i = 1, \dots, \kappa, j = r-\kappa+1, \dots, r}$
be $\kappa \times \kappa$-matrices of indeterminates.

Then
 $$ U \cong \Spec O[A,B]/(AB = BA = \pi) \times V, $$

where 
\begin{eqnarray*}
 V & = & \Spec O[a^0_{i,j}; i = 1, \dots r, j = \kappa+1, \dots, n-r] \times \\
 &&  \Spec O[a^\kappa_{i,j}; i=1, \dots, r-\kappa, j = 1, \dots, \kappa] \\
 & \cong & \A^{(n-r)r-\kappa^2}.
\end{eqnarray*}

{\em Second case: $\kappa > r$}

Let $A=(a^0_{i,j})_{i=\kappa-r+1, \dots, \kappa, j = 1, \dots, r}$,
 $B=(a^\kappa_{i,j})_{i=n-r-\kappa+1, \dots, n-\kappa, j = 1, \dots, r}$
be $r \times r$-matrices of indeterminates.

Then
 $$ U \cong \Spec O[A,B]/(AB = BA = \pi) \times V, $$

where 
\begin{eqnarray*}
  V & = & \Spec O[a^0_{i,j}; i = 1, \dots r, j = \kappa+1, \dots, n-r] \times \\
&&  \Spec O[a^\kappa_{i,j}; i=1, \dots, r-\kappa, j = 1, \dots, \kappa] \\
& \cong & \A^{(n-r)r -r^2}. 
\end{eqnarray*}

So we see that up to a product
with an affine space the special fibre of $U$ is 
 a 'variety of circular complexes'
(cf. \cite{MT}). These varieties have been first studied by Strickland
\cite{Strickland}
using the technique of algebras with straightening law.
She gives an 
explicit $k$-basis
in terms of Young tableaux, and shows that these rings
are algebras with straightening law.
In particular we get 

\begin{thm} \label{strickland}
The special fibre of $\Mloc_{\mu,\kappa}$ is reduced. \qed
\end{thm}

This is the only result we will need, but in fact Strickland proves
much more, for example that the irreducible components of
these varieties are normal and have Cohen-Macaulay singularities.

Recently, Mehta and Trivedi have proved similar
results (for $\Char k = p > 0$) 
using the technique of Frobenius splittings (cf. \cite{MT}).
In fact, their results yield that $\Mlocs_{\mu,\kappa}$
is Frobenius split. But they do not consider $\Mlocs_{\mu, \kappa}$
as a subvariety of the affine flag variety (which is important
for us since we want to look at all the $\Mlocs_{0,\kappa}$ at the
same time and at their intersections). 

\section{Flatness of the Standard Local Model}

\begin{thm} \label{reduced_charp}
Let $\Char k = p > 0$.
Then the special fibre of $\Mloc$ is reduced.
\end{thm}

{\em Proof.}
For $I \subseteq \{ 0,\dots,n-1 \}$, we have the parahoric local model
 $\Mloc_I$. Its special fibre $\Mlocs_I$ can be embedded in
 $SL_n(k\xT)/P^I$, where $P^I$ is a certain
parahoric subgroup of $SL_n(k\xT)$ (cf. section 
\ref{locmod_afffm}).
Denote by $\Mloct_I$
the inverse image of $\Mlocs_I$ under the canonical projection
 $\cF \lto SL_n(k\xT)/P^I$. Set-theoretically,
it is a union of Schubert varieties.

We can describe $\Mloct_I$ in terms of lattice chains in
the following way:
As in section \ref{locmod_afffm}, write
\begin{eqnarray*} 
 &&  \lambda_0 = R\xt^{n}, 
   \lambda_1 = t^{-1} R\xt^{1} \oplus R\xt^{n-1}, \\
 &&\qquad  \dots, 
   \lambda_{n-1} = t^{-1}R\xt^{n-1} \oplus R\xt^{1},
\end{eqnarray*}
where $R$ is a $k$-algebra.
Then $\Mlocs(R)$ is identified via the embedding $\ii$ with 
the set of those complete lattice chains $(\sL_i)_i \in \cF(R)$
with $\lambda_i \supseteq \sL_i \supseteq t\lambda_i$
for all $i$. The closed subscheme $\Mloct_I \subseteq \cF$ consists of the
complete lattice
chains $(\sL_i)_i$ such that 
 $\lambda_i \supseteq \sL_i \supseteq t\lambda_i$ holds for $i \in I$.

Of course, $\Mloc_\kappa$, $\kappa \in \{0,\dots,n-1\}$, 
is just a Grassmannian,
so in particular a smooth variety. Hence $\Mloct_\kappa$
is smooth as well by proposition \ref{projsmooth}. But it is also 
connected, hence irreducible, and hence it is
just a single, smooth Schubert variety. 

Furthermore, the $\Mloct_{\mu,\kappa}$ are reduced by theorem
\ref{strickland} and proposition \ref{projsmooth}.
Thus they are unions of Schubert varieties also scheme-theoretically.

Obviously, we have (scheme-theoretic intersection inside $\Mloct_0$):
 $$ \Mlocs = \bigcap_{\kappa = 1}^{n-1} \Mloct_{0,\kappa}. $$

Now we apply corollary \ref{schubvar_fsplit} to $\Mloct_0$.
We get that $\Mloct_0$ is Frobenius split and that all the 
 $\Mloct_{0,\kappa}$ are simultaneously compatibly split.
Thus their intersection is reduced by lemma \ref{split_inters}
and proposition \ref{splitred}. \qed

{\bf Remark.} The results of Strickland that we used to prove the theorem
are not really in the spirit of this paper. One could hope that
the theorem can also be proved in the following, more elegant way.

What we have used is that $\Mlocs$ is the intersection
of the $\Mloct_{0,\kappa}$ inside $\Mloct_0$. 
But of course, 
 $\Mlocs$ is also the intersection of the $\Mloct_\kappa$ 
(inside $\cF$, say).

Now if we could find a Schubert variety $X$ in $\cF$ which contains
all the $\Mloct_\kappa$, and which we knew to be normal, we 
could apply the same reasoning as above. Not only would this
be independent of Strickland's results, 
it would even imply the results of Strickland (at least the
part we make use of).

One way to get such an $X$ would be to find a 'sufficiently big'
normal Schubert variety in the affine Grassmannian $SL_n/P^0$,
and to take the inverse image in $\cF$. But even in the Grassmannian,
it seems difficult to see which Schubert varieties are normal.

By a constructibility argument, we can show that the
previous theorem holds also in characteristic 0.

\begin{thm} Let $k$ be of characteristic 0. Then the
special fibre of $\Mloc$ is reduced.
\end{thm}

{\em Proof.}
We first construct a $\Z$-scheme the fibre of which over 
a prime $p$
is the special fibre of the local model in characteristic $p$.

Let $\Lambda_i = \Z^n$, $i= 0, \dots, n-1$ 
with basis $e^i_1, \dots, e^i_n$
and define maps 
$$\phi_i : \Lambda_i \lto \Lambda_{i+1},\quad  
 e^i_j \mapsto \left\{
\begin{array}{ll} e^{i+1}_j & j \ne i+1\\
                     0      & j= i+1
\end{array} \right. $$
($i=0, \dots, n-1$, $\Lambda_n := \Lambda_0$).
Denote by $M$ the corresponding scheme of compatibly chosen
subspaces of rank $r$.

Then the geometric fibre of $M$ over a prime $p$ is just the 
special fibre of the local model in characteristic $p$,
and the geometric generic fibre of $M$ is the special fibre
of the local model in characteristic 0.
Denote by $f$ the morphism $M \lto \Spec \Z$.

We know already that all the geometric fibres over primes $p$ are reduced.
As $M$ is a Jacobson scheme, the union of all these
fibres is very dense in $M$. Thus the only constructible
subset of $M$ that contains all these fibres is $M$ itself.
Now the proposition follows from \cite{EGAIV} 9.9.2,
which asserts that in our situation the set of $x \in M$
such that $f^{-1}(f(x))$ is geometrically reduced in $x$
is constructible.
\qed

As a consequence of the previous theorems we get for $k$ of
arbitrary characteristic:

\begin{thm} \label{flatness}
The standard local model $\Mloc$ is flat over $O$.
\end{thm}

{\em Proof}.
The irreducible components of the special fibre correspond 
to the extreme alcoves $x \in W(\mu)$ (see prop. \ref{propsUx}). 
A non-empty open subset of
the irreducible component
associated to $x$ is contained in the open subset $U_x \subseteq \Mloc$,
which is isomorphic to affine space 
as we have seen in proposition \ref{Ux_smooth}.
Thus its generic point can be lifted to the generic fibre.
Now the flatness is a consequence of the reducedness
of the special fibre.
\qed

More generally,

\begin{kor} All the parahoric local models
 $\Mloc_I$, $I \subseteq \{1,\dots,n\}$, are flat over $O$.
\end{kor}

{\em Proof.}
The scheme $\Mloct_I$ is smooth
over $\Mlocs_I$, as it is just the inverse image 
of $\Mlocs_I$ under the morphism $\cF \lto SL_n/P^I$
which is smooth.

As $\Mloct_I$ is reduced by what we have seen
above, so is $\Mlocs_I$.
Now it follows that $\Mloc_I$ is flat over $O$ in the same way
as above. \qed

\subsection{The Singularities of the Local Model}

From the corresponding results about affine Schubert varieties
(see section \ref{sing_schub_var}) we obtain

\begin{stz} \label{singularities}
The irreducible components of $\Mlocs$ are normal. They
have rational singularities, hence in particular are Cohen-Macaulay.
\qed
\end{stz}

It would be very interesting to know if the whole special
fibre of the local model is still Cohen-Macaulay.  

{\bf Remark.} As computations with the help of
computer algebra programs show,
in general the irreducible components
of the special fibre are not locally complete intersections.

\section{The General Case}
\label{unram_ext}

In this section we want to show that our results carry over
without any difficulties to the case where an unramified
field extension is involved in the (EL) datum. We use the same
notation as in \cite{RZ}.
First we define the notion of unramified data of (EL) type;
compare the introduction.

So consider a finite {\em unramified} extension $F/\Q_p$,
and let $B=F$, $V=F^n$. The algebraic group associated
to these data is $G = \Res_{F/\Q_p} GL_F(V)$.

Furthermore, choose an algebraically closed field $L$ 
of characteristic $p$,
and let $K_0$ be the quotient field of the Witt ring $W(L)$.
Since the extension $F/\Q_p$ is unramified, with the notation 
of \cite{RZ} we have $K_0=K$.
 
Now let $\mu: \G_{m,K} \lto G_K$ be a 1-parameter subgroup,
such that the weight decomposition of $V \tens_{\Q_p} K$ contains
only the weights 0 and 1: 
$$V\tens_{\Q_p} K = V_0 \oplus V_1.$$
We denote by $E$ the field of definition of the conjugacy
class of $\mu$.

Finally, let $\sL$ be a periodic lattice chain of lattices 
in $V$.

We will call data of this type unramified data of (EL) type. 
This is a more general notion than is used in 
\cite{RZ}, 3.82. There, the lattice chain consists only of
multiples of 
one lattice and the resulting local model is smooth, whereas
here 'unramified' relates only to the field extension.

The local model $\Mloc$ associated to these
data is the $O_E$-scheme defined as follows. For an $O_E$-scheme $S$,
the $S$-valued points are given by 
\begin{enumerate}
\item a functor $\Lambda \mapsto t_\Lambda$ from $\sL$
      to the category of $O_F \tens_{\Z_p} \cO_S$-modules 
      on $S$, and
\item a morphism of functors $\varphi_\Lambda : \Lambda \tens_{\Z_p}
      \cO_S \lto t_\Lambda$,
\end{enumerate}
which are subject to the following conditions

a) $t_\Lambda$ is locally on $S$ a free $\cO_S$-module of finite
   rank, and we have an identity of polynomial functions on $O_F$:
\begin{equation} \label{kottwitz_cond}
 \det_{\cO_S}(a; t_\Lambda) = \det_K(a;V_0).
\end{equation}

b) $\varphi_\Lambda$ is surjective for all $\Lambda$. 

We have a decomposition
$$ F \tens_{\Q_p} K = \bigoplus_{\varphi: F \lto K} K, $$
and correspondingly we get
$$ V \tens_{\Q_p} K = \bigoplus_\varphi V_\varphi, \quad
   V_0 = \bigoplus_\varphi V_{0,\varphi}, \quad
   V_1 = \bigoplus_\varphi V_{1,\varphi}. $$

All the $V_\varphi$ have dimension $n = \dim_F V$ over $K$.
We write $r_\varphi = \dim V_{0,\varphi}$.
Of course, the number of summands is $d = [F:\Q_p]$.

As $F/\Q_p$ is unramified, we even have a ring
isomorphism
$$ O_F \tens_{\Z_p} O_K = \bigoplus_\varphi O_K. $$

So a $O_F \tens_{\Z_p} O_K$-module $M$
is just a family $(M_\varphi)_\varphi$ of 
 $O_K$-modules, and homomorphisms $M \lto N$ 
are families $(M_\varphi \lto N_\varphi)_\varphi$
of homomorphisms.

To investigate properties like flatness or reducedness
of the special fibre, we can just as well look at 
 $\Mloc \tens_{O_E} O_K$.

We get the following description: For an $O_K$-scheme $S$,
the $S$-valued points of the local model are given by

\begin{enumerate}
\item a functor which associates to each $\Lambda$ in $\sL$
      a family $(t_{\Lambda_\varphi})$ of $\cO_S$-modules, and
\item a morphism of functors 
$$\varphi_\Lambda : \Lambda \tens_{\Z_p}
      \cO_S = \bigoplus_\varphi \Lambda_\varphi \tens_{O_K} \cO_S 
      \lto \bigoplus_\varphi t_{\Lambda,\varphi},$$
\end{enumerate}
which are subject to the following conditions

a) $t_{\Lambda,\varphi}$ is locally on $S$ a free $\cO_S$-module of finite
   rank $r_\varphi$, and

b) all the morphisms $\Lambda_\varphi \tens_{O_K} \cO_S 
      \lto t_{\Lambda,\varphi}$ are surjective.

Thus the local model in this case is just a product of 
standard local models. In particular, we get
from theorem \ref{flatness} and proposition \ref{singularities} 
our main theorem.

\begin{thm} \label{general_case}
The local model $M$ associated to an unramified (EL)-datum 
is flat over $O_E$, and its special fibre is reduced.
Furthermore, the irreducible components of the special fibre
of $M$
are normal with rational singularities,
so in particular are Cohen-Macaulay. \qed
\end{thm}

\vskip2cm
\hspace*{1cm}Ulrich G\"ortz\\
\hspace*{1cm}Mathematisches Institut\\
\hspace*{1cm}der Universit\"at zu K\"oln\\
\hspace*{1cm}Weyertal 86--90

\hspace*{1cm}{\small DE--}50931 K\"oln (Germany)

\hspace*{1cm}ugoertz@mi.uni-koeln.de

\end{document}